\documentclass[12pt]{article}

\usepackage{latexsym,amsmath,amscd,amssymb,graphics}
\usepackage{diagrams}
\newarrow {Corresponds} <{dash}{}{dash}>
\usepackage{enumerate}
\usepackage{framed}
\usepackage{graphicx}
\usepackage{eucal}
\usepackage[square,authoryear]{natbib}
\usepackage[colorlinks]{hyperref}
\usepackage{url}
\usepackage{mathrsfs}
\usepackage{mathabx}
\usepackage{ulsy}

\usepackage[all]{xy}

\makeatletter

\@addtoreset{figure}{section}
\def\thefigure{\thesection.\@arabic\c@figure}
\def\fps@figure{h, t}
\@addtoreset{table}{bsection}
\def\thetable{\thesection.\@arabic\c@table}
\def\fps@table{h, t}
\@addtoreset{equation}{section}

\makeatother

\begin{document}

\newtheorem{theorem}{Theorem}[section]
\newtheorem{definition}[theorem]{Definition}
\newtheorem{lemma}[theorem]{Lemma}
\newtheorem{remark}[theorem]{Remark}
\newtheorem{proposition}[theorem]{Proposition}
\newtheorem{corollary}[theorem]{Corollary}
\newtheorem{example}[theorem]{Example}

\def\below#1#2{\mathrel{\mathop{#1}\limits_{#2}}}



\title{Lagrangian reductions and integrable systems in 
condensed matter}
\author{Fran\c{c}ois Gay-Balmaz$^{1}$, Michael Monastyrsky$^{2}$, Tudor S. Ratiu$^{3}$}
\addtocounter{footnote}{1} 
\footnotetext{Laboratoire de 
M\'et\'eorologie Dynamique, \'Ecole Normale Sup\'erieure/CNRS, Paris, France. Partially supported by ``Projet incitatif de recherch de l'ENS", by the government grant of the Russian Federation for support of research projects implemented by leading scientists, Lomonosov Moscow State University under the agreement No. 11.G34.31.0054, and the 
Isaac Newton Institute for Mathematical Sciences, Cambridge, UK.
\texttt{gaybalma@lmd.ens.fr}
\addtocounter{footnote}{1} }
\footnotetext{Institute of Theoretical and Experimental Physics, Moscow, Russia. Partially supported by the government grant of the Russian Federation for support of research projects implemented by leading scientists, Lomonosov Moscow State University under the agreement No. 11.G34.31.0054, RFBR grant 13-01-00314, and the Isaac Newton Institute for Mathematical Sciences, Cambridge, UK.
\texttt{monastyrsky@itep.ru}
\addtocounter{footnote}{1} }
\footnotetext{Section de Math\'ematiques and Bernoulli Center, 
\'Ecole Polytechnique F\'ed\'erale de
Lausanne. CH--1015 Lausanne. Switzerland. Partially supported by Swiss NSF grant 200020-126630, by the government grant of the Russian Federation for support of research projects implemented by leading scientists, Lomonosov Moscow State University under the agreement No. 11.G34.31.0054, and the 
Isaac Newton Institute for Mathematical Sciences, Cambridge, UK.
\texttt{tudor.ratiu@epfl.ch}
\addtocounter{footnote}{1} }

\date{}
\maketitle

\makeatother
\maketitle




\begin{abstract}
We consider a general approach for the process of Lagrangian and Hamiltonian reduction by symmetries in 
chiral gauge models. This approach is used to show the complete integrability of several one dimensional texture equations arising in liquid Helium phases and neutron stars.
\end{abstract}

\tableofcontents


\section{Introduction}

There is a well established relation between quantum field theory and
condensed matter physics. For example, physical phenomena such a superfluidity
and superconductivity are the manifestation of quantum effects at
microscopic level. On the other hand, there are classical systems,
such as liquid crystals, where many phenomena, such as phase transitions
between different mesophases, are described in the framework of the
Landau-de Gennes theory. These systems include superfluid ${}^3$He, 
the superfluid core of neutron stars, biaxial and uniaxial nematics.

The common feature of these different systems is that in some interval 
of the transition temperature, their behavior is determined by the 
Ginzburg-Landau equation with multidimensional order parameters. Another
interesting feature of these systems is the existence of different
thermodynamic phases. The description of phase transitions between 
different phases is a difficult and important problem in 
condensed matter physics. The approach, based on the identification of
thermodynamic phases with orbits of the group of symmetry of the 
potential in the free energy, as developed in \cite{GoMoNo1979}, 
\cite{GoMo1978a, GoMo1978b}, and \cite{BoMo1987}, is a useful tool 
for a complete classification of phases and gives a global description 
of these phases. From many points of view, it is possible to study 
these systems like a chiral field model. Viewed this way, it is
known that these systems can be obtained by a general reduction procedure 
developed in \cite{GBTr2010}, based on \cite{CeMaRa2001}, \cite{HoMaRa1998}, \cite{CeMaPeRa2003}. The goal of this paper is to unify these two approaches
and techniques, to formulate a theory that contains them both, and
especially to show its effectiveness by studying in detail the
complete integrability
of several concrete physical systems in different phases.

We begin with a short review of the relevant facts of the 
Lagrange-Poincar\'e and Euler-Poincar\'e variational principles in 
Section \ref{sec:general_theory}. Its Hamiltonian counterpart, 
Hamilton-Poincar\'e and Lie-Poisson reduction, are treated in
Section \ref{sec:HPandLP}. We shall limit ourselves with the classical, as opposed to
the field theoretical, description of these theories, because all
examples analyzed in this paper necessitate only this classical 
theory. The field theoretical approach, which we have also developed,
will be the subject of another paper.
The main result of these sections is an equivalence of the two 
descriptions. These results are used in Section \ref{sec:applications},
forming the mani body of the paper, to study in detail the behavior 
of superfluid ${}^3$He and neutron star cores in different phases. 
The versatility of passing from one description to another, 
as well as between the Lagrangian and Hamiltonian formulations,
is crucial in the proof of the complete
integrability of the equations associated to different phases. The
key to the success of our geometric method is the fact that all 
physical systems under study have a natural Lagrangian and 
Hamiltonian formulation within the Lagrange-Poincar\'e and 
Hamilton-Poincar\'e theories, with the Lagrangian and Hamiltonian
independent on a very special group of variables. This implies
that these systems have an equivalent Euler-Poincar\'e and
Lie-Poisson description which turns out to be considerably
simpler and more appropriate to the study of the dynamics of the
equations associated to the relevant phases. The possibility of 
using at once the four descriptions of the systems under consideration
leads directly to the proof of complete integrability of the
equations describing the system's behavior in different phases.


\section{Lagrange-Poincar\'e and  Euler-Poincar\'e reduction on Lie groups}
\label{sec:general_theory}
 
In this section we shall quickly review two Lagrangian reduction processes, namely
Lagrange-Poincar\'e and Euler-Poincar\'e reduction, as they apply to a Lagrangian defined on a Lie group and invariant under right translation by a closed subgroup. We shall also emphasize the case of discrete
symmetry groups.

\paragraph{Geometric setup.} Let $M$ be the parameter manifold of 
the theory and let $\Phi: G \times M \rightarrow M$ be a left 
transitive Lie group action. Usually, $M$ is a particular 
orbit of the action of $G$ on a bigger manifold. Selecting one 
particular orbit corresponds to choosing a particular phase of 
the physical system. Given a Lie group $G$, we shall denote
by the corresponding Fraktur letter $\mathfrak{g}$ is Lie algebra.

Choose an element $m _0 \in M$ and consider the isotropy 
subgroup $H:=G_{m_0}$. We have the diffeomorphism 
$G/H\ni [g]:=gH \stackrel{\sim}\longmapsto gm_0 \in M$, where 
$H$ acts on $G$ by right multiplication $R_hg := gh$
for all $h \in H$ and $g \in G$. We shall always identify 
$M$ with $G/H$ via this diffeomorphism 
and denote by $\pi:G 
\rightarrow G/H$ the orbit space projection.

We suppose that the theory is described by a Lagrangian
$\mathfrak{L} = \mathfrak{L} (m, \dot m):TM \rightarrow \mathbb{R}$, whose associate Euler-Lagrange equations read
\[
\frac{d}{dt} \frac{\partial \mathfrak{L}  }{\partial \dot m}- \frac{\partial \mathfrak{L} }{\partial m}=0.
\]
Recall that these equations follow from applying Hamilton's principle
\[
\delta \int_{ t _0 }^{ t _1} \mathfrak{L}  ( m(t), \dot m(t))dt=0,
\]
for arbitrary variations of the curve $m (t)$ whose corresponding
infinitesimal variations $\delta m(t)$ satisfies
$ \delta m(t _0 )= \delta m( t _1 )=0$.

Since the $G$-action is transitive on $M = G/H$, any curve 
$m: [ t _0 , t _1 ]  \rightarrow M$ can be written as $m(t)= 
\Phi _{g(t)} (m_0) = : g(t) m_0$, where $g: [ t _0 , t _1 ] \rightarrow G$.
By using this relation we can rewrite the action functional in terms of the curve $g (t) $ as
\begin{align*} 
\int_{ t _0 }^{ t _1 } \mathfrak{L}  (m(t) , \dot m(t) )dt &=\int_{ t _0 }^{ t _1 } \mathfrak{L} \left( g (t) m _0 , \frac{d}{dt} g(t) m _0 \right) dt \\
&= \int_{ t _0 }^{ t _1 } \mathfrak{L}  \left( g (t) m _0 ,( \dot g(t)g(t) ^{-1} )_M( g(t) m _0 )\right) dt,
\end{align*}
where, for every $\xi\in \mathfrak{g}$, 
$M\ni m \mapsto \xi_M(m): = 
\left.\frac{d}{dt}\right|_{t=0}\Phi_{\exp t \xi}(m) \in T_mM$ denotes
the infinitesimal generator vector field of the action.
This suggests the definition of the Lagrangian $L_{m_0}$ for curves $g(t)$ in the Lie group as
\[
L_{m _0 } :TG \rightarrow \mathbb{R} , \quad 
L _{ m _0}( g, \dot g) := \mathfrak{L} 
\left( g m _0 ,( \dot gg ^{-1} )_M( g m _0)\right).
\]
This Lagrangian is clearly $H$-invariant. Our goal is to find an explicit relation between the 
Euler-Lagrange equations for $L_{m _0 }$ and $ 
\mathfrak{L} $ as well as to deduce another simpler equivalent 
formulation of these equations.

\medskip

To do this, we shall start with a $H$-invariant Lagrangian 
$L :TG \rightarrow \mathbb{R}  $ and, following \cite{GBTr2010}, 
we shall carry out two 
reductions processes for $L$. The first one follows the 
Lagrange-Poincar\'e reduction theory (see \cite{CeMaRa2001}), 
the second one is a generalization of the
Euler-Poincar\'e reduction with parameters (see \cite{HoMaRa1998}). 
These two reductions  
correspond to two realizations of the quotient space $(TG)/H$.

\paragraph{Lagrange-Poincar\'e approach.} The 
Lagrange-Poincar\'e reduction is implemented by using the vector bundle 
isomorphism
\[
\alpha _\mathcal{A} : (TG)/H \rightarrow TM
\times_{M} \widetilde{\mathfrak{h}}
\]
over $M=G/H$. Here
$\widetilde{\mathfrak{h}} := G \times _H \mathfrak{h} 
\rightarrow M$ is the adjoint bundle, where the
right $H$-action on $G \times \mathfrak{h}$ is given by 
$(g, \eta)\cdot h: = (gh, \operatorname{Ad}_{h^{-1}} \eta)$
for all $h \in H$, $g \in G$, and $\eta\in \mathfrak{h}$.
The vector bundle isomorphism $\alpha _{\mathcal{A}}$ is constructed with the help of a principal connection $\mathcal{A}\in 
\Omega^1 (G, \mathfrak{h})$ on the principal bundle 
$\pi : G \rightarrow  G/H =M$ and reads
\[
\alpha _\mathcal{A} ([v _g ]_H):= \left(T_{g}\pi (v _g ),
\left[g, \mathcal{A} (v _g )\right]_ H \right)=\left( (v _g g^{-1} )_M(m),\left[g, \mathcal{A} (v _g )\right]_ H \right)
\]
(we denote by $[x]_H$ a point in orbit space of the $H$-action on the manifold whose points are $x$).

From the given right $H$-invariant Lagrangian 
$L:TG \rightarrow \mathbb{R}$, 
we get the reduced Lagrangian $\mathcal{L} :TM \times _M \widetilde{\mathfrak{h}} \rightarrow \mathbb{R}  $,  $\mathcal{L} =\mathcal{L} (m, \dot m, \sigma )$, defined by
\begin{equation}
\label{LP_Lagrangian}
L(g, \dot g)=\mathcal{L} \left(  gm _0 , \dot g m _0 ,\left[g, \mathcal{A} (g, \dot g)\right] _H \right) .
\end{equation}

The reduced Euler-Lagrange equations (or Lagrange-Poincar\'e equations) are obtained by computing the critical curve of the variational principle
\begin{equation}\label{LP_VP} 
\delta \int_{ t _0 }^{ t _1 } \mathcal{L} (m, \dot m, \sigma )dt,
\end{equation} 
for variations $ \delta m$ and $ \delta ^ \mathcal{A} \sigma $ induced by variations $ \delta g(t)$ of the curve $g(t) $, that vanish at $t=t_0,t_1$. While the variations $ \delta m(t)$ are free and vanish at $t=t _0 , t _1 $, the variations of $ \sigma (t) $ verify
\begin{align*} 
\delta^\mathcal{A} \sigma (t)   &:= 
\left.\frac{D^\mathcal{A} }{D\varepsilon}\right|_{\varepsilon=0} 
\left[g_\varepsilon (t), \mathcal{A} (g_\varepsilon (t), \dot g_\varepsilon (t) ) \right]_H \\
&= \frac{D^ \mathcal{A} }{Dt} \eta (t) +[ \eta (t) , \sigma (t)   ]+  \mathbf{i} _{ \delta m (t) } \widetilde{ \mathcal{B} } \in \widetilde{\mathfrak{h}},
\end{align*} 
where $D^\mathcal{A}/D\varepsilon$ denotes 
the covariant derivative defined by the connection one-form 
$\mathcal{A}$,
$\widetilde{\mathcal{B}} \in \Omega^2 (M, \widetilde{\mathfrak{h}})$ 
is the reduced curvature on the base associated to the connection 
$ \mathcal{A}$, and   
$\eta (t)= \left[g (t), \mathcal{A} (\delta g(t)) \right]_H \in 
\widetilde{\mathfrak{h}}$ is arbitrary with $ \eta (t_0)= \eta (t _1 )=0$. Using these variations in \eqref{LP_VP} yield the
Lagrange-Poincar\'e equations
\begin{equation}\label{LP_equations} 
\frac{D^ \mathcal{A} }{Dt} \frac{\delta \mathcal{L} }{\delta \sigma }+ \operatorname{ad}^*_ \sigma   \frac{\delta \mathcal{L} }{\delta \sigma }=0, \quad\quad  \frac{\partial  \mathcal{L} }{\partial  m}- \frac{d}{dt} \frac{\partial  \mathcal{L} }{\partial  \dot m}= \left\langle \frac{\delta \mathcal{L} }{\delta \sigma }, \mathbf{i} _{\dot m} \widetilde{ \mathcal{B} } \right\rangle.
\end{equation}

We refer to \cite{CeMaRa2001}  for the general theory and to \cite{GBTr2010} for this special case.

\paragraph{Lagrange-Poincar\'e equations for $H$ discrete.} Assume now that
$H$ is a closed discrete subgroup of $G$. Then $\mathfrak{h} = \{0\}$,
$\widetilde{\mathfrak{h}}$ is the vector bundle with zero dimensional
fiber and base $M$, $\mathcal{A} = 0$, and hence the vector bundle
isomorphism $\alpha _\mathcal{A}$ becomes canonical, $\alpha : (TG)/H \rightarrow TM$, the source and target spaces viewed as a vector bundles over $M$, and it is given
by
\begin{equation*} 
\alpha ([v _g ]_H):= T _g\pi ( v _g )=
(v _g g ^{-1} )_M(m)\in T_{m}M.
\end{equation*}
So, if $g:[ t _0 , t _1 ] \rightarrow G$ is a given curve, we get simply
\begin{equation*} 
\alpha\left([g(t), \dot g(t)]_H \right) =\left( g(t) m _0 , \frac{d}{dt} g(t) m _0 \right)  = 
(m(t), \dot m(t))\in TM.
\end{equation*}
The reduced Lagrangian $\mathcal{L} :TM\rightarrow \mathbb{R}$
yields the Lagrange-Poincar\'e equations \eqref{LP_equations} which in
this case  become
\begin{equation}
\label{LP_equations_discrete}
\frac{\partial\mathcal{L} }{\partial  m}- \frac{d}{dt}  
\frac{\partial\mathcal{L} }{\partial \dot m }=0.
\end{equation} 

It is instructive to consider in more 
detail the isomorphism $\alpha $ in the case of a closed 
discrete subgroup. In this case, the kernel of the tangent 
map is zero, so that at any $g \in G$ we have the isomorphism 
$T_g \pi :T_g G \rightarrow T_{[g]}(G/H)$ which implies that
\[
\alpha : (TG)/H \rightarrow T(G/H), \quad [v_g]_H \mapsto 
T \pi (v_g)
\]
is a vector bundle isomorphism covering the identity on $G/H$. 
Indeed, if $T_g \pi (v _g )= T _{\bar g}\pi ( w _{\bar g})$, 
then necessarily $\bar g =gh$ for $ h \in H$, and we can write 
$T _g\pi (v _g )= T_{\bar g} \pi (w_{\bar g})= 
T_g \pi (w_{\bar g} h^{-1})$, so that $v_g = 
w_{\bar g} h ^{-1}$, since $\operatorname{ker} T_g\pi =\{0\}$. 
This proves that $[v_g ]_H = [w_{\bar g} ]_H$.

\paragraph{Euler-Poincar\'e approach.} The 
Euler-Poincar\'e reduction is implemented by using the vector bundle 
isomorphism
\[
\bar i_{m_0 }: (TG)/H \rightarrow \mathfrak{g}\times M, \quad \bar i_{m_0 }( [v _g ]_H)= \left( v _g g ^{-1} , \Phi _g (m_0 )\right) ,
\]
where an element $m_0 \in M$ has been fixed, \cite{GBTr2010}.
We note that a connection is not needed to write this isomorphism. 
If $g:[ t _0 , t _1 ] \rightarrow G$ is a given curve, 
this formula implies
\begin{equation}
\label{reduced_section_formula_for_g}
\bar i_{m_0 } \left([g(t) , \dot g(t) ]_H \right)= 
\left( \dot g (t) g (t) ^{-1} , \Phi _{g(t)}(m_0)\right) = :
(\xi (t), m(t)).
\end{equation}

Note that by composing the two vector bundle isomorphisms 
$\alpha _{\mathcal{A}}$ and $\bar i_{m_0}$ over $M$, 
we get the vector bundle isomorphism
\begin{equation}\label{isomorphism_long} 
\mathfrak{g}\times M\ni (\xi,m) \longmapsto 
\left(\xi _M(m), \left[g,\mathcal{A}(\xi g)\right]_H\right) 
\in  TM  \times \widetilde{ \mathfrak{h}},
\end{equation} 
over $M$, where $g \in G$ is arbitrary such that $\pi (g)=m$.

Given $v _m \in  T_mM$, $\xi _m \in  \tilde{\mathfrak{h}} _m$, the inverse of the above map is given by
\begin{equation}\label{isomorphism_inverse_long} 
TM \times \widetilde{ \mathfrak{h}} \ni (v _m , \xi _m )  \longmapsto \Big( \left( \operatorname{Hor}_g(v _m ) \right) g ^{-1} + \operatorname{Ad}_g \eta, m \Big)    \in \mathfrak{g} \times M,
\end{equation}  
where $g \in G$ is such that $ \pi (g)=m$, $ \operatorname{Hor}_g:T_mM \rightarrow T_gG$ is the horizontal lift of the connection $ \mathcal{A} $, and   $ \eta \in  \mathfrak{h} $ is such that $\xi _m =[g, \eta ]_H$.
A direct verification shows that this expression does not depend on $g$ as long as $ \pi (g)=m$.

Given a $H$-invariant Lagrangian $L:TG \rightarrow \mathbb{R}$, the associated 
reduced Lagrangian $l: \mathfrak{g}\times M \rightarrow 
\mathbb{R}$ obtained through the Euler-Poincar\'e process is
\begin{equation}
\label{EP_Lagrangian}
L(g, \dot g)=l\left(\dot g g^{-1}, \Phi_g(m_0) \right)=l( \xi , m).
\end{equation}
The Euler-Poincar\'e equations for $\mathfrak{L}$ follow from applying the variational principle with constrained variations
\[
\delta \int_{ t _0 }^{ t _1 }l( \xi , m) dt=0, \quad 
\delta \xi = \dot{\eta} +[ \eta  , \xi  ] , \quad \delta m= \eta  _M(m),
\]
where $ \eta(t) \in  \mathfrak{g}$ is arbitrary curve with $ \eta ( t _0 )= \eta ( t _1 )=0$. We thus get the equations
\begin{equation}
\label{EP_equations}
\frac{d}{dt} \frac{\delta l}{\delta \xi }+\operatorname{ad}^*_\xi  \frac{\delta l}{\delta \xi }= \mathbf{J} \left(  \frac{\delta l}{\delta m } \right), \quad  \dot m= \xi _M (m),
\end{equation}
where $\mathbf{J}: T^\ast M \rightarrow \mathfrak{g}^\ast$ is the standard equivariant momentum map of the cotangent lifted action
given by $\left\langle\mathbf{J}(\alpha_m), \zeta \right\rangle = 
\left\langle\alpha_m, \zeta_M(m)\right\rangle$ for all $\alpha_m\in 
T^\ast_mM$, $\zeta\in \mathfrak{g}$. 

Note that if $H$ is closed and discrete, the Euler-Poincar\'e equations \eqref{EP_equations} for $l:\mathfrak{g} \times M \rightarrow \mathbb{R}$ do not simplify, contrary to what happens on the Lagrange-Poincar\'e side.

\section{Hamilton-Poincar\'e and Lie-Poisson reduction on Lie groups}\label{sec:HPandLP} 

In this section we will summarize the necessary material that is relevant for the Hamiltonian description of condensed matter systems.
This is the Hamiltonian version of the two Lagrangian reduction processes described in the preceding section.

\paragraph{Hamilton-Poincar\'e reduction.} As in the preceding section, we consider a Lie group acting transitively on the left on a manifold $M$. Choosing $m_0 \in M$, we have the diffeomorphism $G/H\ni gH \mapsto \Phi _g ( m _0 ) \in M$. Given a Hamiltonian $H : T^* G \rightarrow \mathbb{R}  $ that is right $H$-invariant, we obtain, by reduction, a Hamiltonian defined on the quotient space $(T^*G)/H$. Similarly as before, choosing a principal connection $ \mathcal{A} \in \Omega ^1 (G, \mathfrak{h}  )$, we have a vector bundle isomorphism
\[
(T^*G)/H \rightarrow T^*(G/H) \oplus \widetilde{\mathfrak{h}} ^\ast , \quad [ \alpha _g ]_H  \mapsto \left( \operatorname{Hor}_g ^\ast \alpha _g , [ g, \mathbf{J} ( \alpha _g )]_H \right) =:( \alpha _m , \bar \mu ),
\]
where $ \operatorname{Hor}_g ^\ast :T_mM \rightarrow T^*G$ is the dual map to the horizontal lift $ \operatorname{Hor}_g:T_mM \rightarrow T_gG$ associated to the connection $ \mathcal{A} $, and $ \mathbf{J} :T^*G \rightarrow \mathfrak{h}  ^\ast $ is the momentum map associated to right translation by $H$. The reduced Hamilton equations obtained by Poisson reduction are called the Hamilton-Poincar\'e equations and read
\begin{equation}\label{Ham_Poincare} 
\frac{D y}{Dt} =- \frac{\partial \mathcal{H} }{\partial x}- \left\langle \bar \mu  , \tilde{ \mathcal{B} }( \dot x, \_\, ) \right\rangle , \quad \dot x= \frac{\partial \mathcal{H} }{\partial y}, \quad \frac{D^ \mathcal{A} \bar \mu }{Dt}  +\operatorname{ad}^*_{ \frac{\delta \mathcal{H} }{\delta \bar \mu }} \bar \mu=0,
\end{equation} 
where $(x,y) \in T^*(G/H)$, $\bar \mu \in \widetilde{\mathfrak{h}} ^\ast $, $\tilde{ \mathcal{B} }\in \Omega ^2 ( G/H, \widetilde{ \mathfrak{h}  })$ is the reduced curvature of $ \mathcal{A} $, and $D/Dt$ in the first equation denotes the covariant derivative on $T^*(G/H)$ associated to a given affine connection on $G/H$, and 
$D^ \mathcal{A}/Dt$ in  the last equation denotes the covariant derivative on $\widetilde{\mathfrak{h}} ^\ast$ associated to the principal connection $ \mathcal{A} $; for details see \cite{CeMaPeRa2003}.

The symplectic leaves in $ T^*(G/H) \oplus \widetilde{\mathfrak{h}} ^\ast $ have been described in \cite{MaPe2000}; they are of the form
$T^*(G/H) \times _{ G/H} \widetilde{ \mathcal{O}}$,
where $ \mathcal{O} $ is a coadjoint orbit of $H$, and 
$ \widetilde{ \mathcal{O} }\rightarrow G/H$ is the associated 
fiber bundle. The symplectic form is the sum of the canonical 
symplectic form on $T^*(G/H)$ and a two-form on 
$\widetilde{ \mathcal{O} }$, see \cite[Theorem 2.3.12]{MaMiOrPeRa2007}. 
If the Lie group $G$ is connected and $\mathcal{O}$ has $N$ elements 
(which is happening in subsequent applications), then the 
fiber bundle $T^*(G/H) \times _{G/H} \widetilde{\mathcal{O}}
\rightarrow G/H$ has $N$ connected components, each one of them 
symplectically diffeomorphic to the canonical phase space $T^*(G/H)$.

\paragraph{Lie-Poisson reduction.} A second realization of $T^*(G/H)$ is given by the diffeomorphism
\[
(T^*G)/H \ni [ \alpha _g ]_H \longmapsto ( \alpha _g g ^{-1} , g m_0 ) \in \mathfrak{g}  ^\ast \times M.
\]
The reduced Hamilton equations on this space read
\begin{equation}\label{Lie_Poisson} 
\dot \mu + \operatorname{ad}^*_{ \frac{\delta h}{\delta \mu }} \mu = \mathbf{J} \left(\frac{\delta h}{\delta m}  \right),\quad \dot m=-\left( \frac{\delta h}{\delta \mu }\right) _M(m), 
\end{equation} 
where $ \mathbf{J} :T^*M \rightarrow \mathfrak{g}  ^\ast $ is the momentum map of the cotangent lift of the left action of $G$ on $M=G/H$ and $h: \mathfrak{g}  ^\ast \times M\rightarrow \mathbb{R}  $ is the reduced Hamiltonian. These equations are Hamiltonian relative to the Poisson bracket
\begin{equation}\label{red_Poisson}
\{f,g\}( \mu , m)=  \left\langle \mu, \left[ \frac{\delta f}{\delta \mu }, \frac{\delta g}{\delta \mu }\right] \right\rangle + \left\langle \mathbf{J} \left(\frac{\delta f}{\delta m}  \right) , \frac{\delta g}{\delta \mu }\right\rangle - \left\langle \mathbf{J} \left(\frac{\delta g}{\delta m}  \right) , \frac{\delta f}{\delta \mu }\right\rangle,
\end{equation}  
see \cite{KrMa1987} (in which a more general situation is considered).

We refer to \cite{GBTr2010}, for further details and examples of application of these two reduction processes.

\medskip

We now suppose that $V$ is a representation space of $G$ and we take 
$M= \operatorname{Orb}(a_0) \subset V^\ast$. The induced Lie algebra
representation $\mathfrak{g} \times V \rightarrow V$ is given by the
infinitesimal operator map and is denoted by $\xi v: = \xi_V(v)$,
for any $\xi \in \mathfrak{g}$ and $v\in V$. We consider the semidirect 
product $S= G \,\circledS\, V$ and its Lie algebra $\mathfrak{s} = 
\mathfrak{g}\,\circledS\,V$. The symplectic leaves in 
$\mathfrak{s}^\ast$ are given by the connected components of the 
coadjoint orbits $\mathcal{O}_{(\mu ,a)}$ of $S$. From the formula 
of the coadjoint action 
\begin{equation}
\label{coadj_semidirect}
\operatorname{Ad}^*_{(g,v)^{-1}}(\mu, a) 
=(\operatorname{Ad}^*_{g^{-1}}\mu + v\diamond ga, ga),
\end{equation}
where $(g,v) \in S$ and $(\mu, a) \in \mathfrak{s}^\ast$, we see 
that the symplectic leaves in $\mathfrak{g}^\ast \times M$ are 
$\mathcal{O}_{(\mu, a_0)}$ endowed with the minus orbit 
symplectic form. The diamond operation  $\diamond: V \times V^\ast
\rightarrow \mathfrak{g}^\ast$ in this formula is defined by 
$\left\langle v \diamond a, \xi \right\rangle : = \left\langle a ,
\xi v \right\rangle$, for any $\xi\in \mathfrak{g}$, where the pairing in the left hand side is
between $\mathfrak{g}^\ast$ and $\mathfrak{g}$, whereas in the
right hand side it is between $V^\ast$ and $V$.

These considerations provide the proof of the following theorem.

\begin{theorem}
\label{symplectic_leaves}  
Given $ \mu \in \mathfrak{g}^\ast $ and $a_0 \in V^\ast$, 
we define $\mu_{a_0}:= \mu |_{\mathfrak{g}_{a_0}} \in 
\mathfrak{g}_{a_0}^\ast$, 
where $ \mathfrak{g}_{ a _0 }=\{ \xi \in \mathfrak{g}  
\mid \xi a _0 =0\}=:\mathfrak{h}$. Let $\mathcal{O}_{\mu_{a_0}}
\subset \mathfrak{g}_{a_0}^\ast$ be the coadjoint orbit of 
$H:=G_{a_0}$ through $\mu_{a_0}$. 
The map
\[
\mathfrak{s}^\ast \supset
\mathcal{O}_{(\mu, a_0)}\ni (\alpha_g g^{-1}, g a_0 ) \longmapsto 
(\operatorname{Hor}_g ^\ast \alpha_g , [g, \mathbf{J}(\alpha_g)]_H) 
\in  T^*(G/H) \times_{G/H} \widetilde{\mathcal{O}}_{\mu_{a_0}}
\]
is a symplectic diffeomorphism.
\end{theorem} 

Note that the Theorem states that if $(\mu-\nu )|_{\mathfrak{g}_{a_0}}
=0$, then $\mathcal{O}_{(\mu, a_0)}=\mathcal{O}_{(\nu, a_0)}$.
This can be verified directly by observing that 
$\operatorname{Ad}^*_{(e,v)^{-1} }(\mu, a_0)=(\mu + v \diamond a_0, 
a_0)$ and the map $ V \ni v \mapsto v \diamond a_0 \in 
\mathfrak{g}_{a_0}^\circ$ (the annihilator of
$\mathfrak{g}_{a_0}$ in $\mathfrak{g}^\ast$) is surjective (which is equivalent to 
$\operatorname{ker}(v \mapsto v \diamond a_0)= \mathfrak{g}_{a_0}$). 

We now write explicitly the operator $\alpha_g \mapsto 
\operatorname{Hor}^*_g \alpha_g $ in the particular case 
when there is an Ad-invariant inner product $\gamma$ on 
$\mathfrak{g}$. We extend $\gamma$ by left invariance to 
a Riemannian metric on $G$. This Riemannian metric, also 
denoted $\gamma$, is right invariant. The principal connection 
on the right $H$-principal bundle $G \rightarrow G/H$ associated 
to $\gamma$ has the expression $\mathcal{A}(v_g):=
\mathbb{P}_{a_0}(g^{-1} v_g )$, where $\mathbb{P}_{a_0 }: 
\mathfrak{g}\rightarrow \mathfrak{g}_{a_0}$ is the 
$\gamma$-orthogonal projection. The horizontal lift associated to $ \mathcal{A} $ reads
\begin{equation}\label{Hor_lift}
\operatorname{Hor}_g:T_mM \rightarrow T_gG, \quad \operatorname{Hor}_g(\xi _M (m) )= g  \, \mathbb{P}_{a _0 }^\perp(\operatorname{Ad}_{g ^{-1} } \xi ),
\end{equation} 
where $ \mathbb{P}  _{a _0 }^\perp: \mathfrak{g}  \rightarrow \mathfrak{g}  _{ a _0 }^\perp$ is the $ \gamma $-orthogonal projection and $m=g a _0 $.
We endow $M=G/H$ with the natural induced Riemannian metric, i.e., 
\begin{equation}\label{induced_metric}
\begin{aligned} 
\gamma _M( \xi _M(m), \eta _M (m) ):&= \gamma \left( 
\operatorname{Hor}_g(\xi _M(m)),\operatorname{Hor}_g(\eta _M(m))\right)\\
& = \gamma \left(\mathbb{P}_{a_0}^\perp(\operatorname{Ad}_{g^{-1}}\xi) , 
\mathbb{P}_{a_0}^\perp(\operatorname{Ad}_{g^{-1}} \eta)\right). 
\end{aligned}
\end{equation} 
Using the Riemannian metrics $\gamma$ and $\gamma_M$, we identify 
$TG$ with $T^*G$ and $TM$ with $T^*M$, respectively. With these identifications, we have
\begin{equation}\label{Hor_Lift} 
\operatorname{Hor}^*_g \alpha _g= ( \alpha _g g ^{-1} )_M(m), \quad \alpha _g \in T^*_gG=T_gG.
\end{equation}
We summarize the maps in this discussion in the following
diagram
\[\hspace{-1.2cm}
\begin{diagram}
&   & T ^\ast G  && 
\\
&\ldTo(2,2)& &\rdTo(2,2)&
\\
\mathfrak{g}^\ast \times M  & & \rTo^{\widetilde{\phantom{mmmm}}}  & & (T^\ast G)/H = T ^\ast (G/H)\oplus\widetilde{\mathfrak{h}}^\ast
\\
\uTo& &  &&\uTo
\\
\mathcal{O}_{(\mu_0, a_0)} & & \rTo^{\widetilde{\phantom{mmmm}}}& &  T ^\ast(G/H) \times _{G/H} \widetilde{\mathcal{O}}_{\mu_{a_0}}.
\end{diagram}
\]

\section{Applications to condensed matter}\label{sec:applications} 

\subsection{Setup of the problem}

\paragraph{Lagrangian description.} As discussed at the beginning of the previous section, for condensed matter theories the Lagrangian $\mathfrak{L} $ is defined on the tangent bundle $TM$ of the parameter manifold $M$. This manifold is assumed to be a homogeneous space, relative to the transitive action of a Lie group $G$, with isotropy group $H=G_{m_0 }$ for some preferred element $m _0 \in M$. From $\mathfrak{L} $ one can construct a Lagrangian $L_{ m _0 }:TG \rightarrow \mathbb{R}$, $L_{ m _0 }(g, \dot g):=\mathfrak{L}  ( g m _0 , (\dot g g ^{-1} )_M( g m _0 ))$ as explained earlier.

Using the results of \S\ref{sec:general_theory}, we will show that the Euler-Lagrange equations for $L_{ m _0 }$ are equivalent to those for 
$\mathfrak{L} $ by implementing Lagrange-Poincar\'e reduction. Then 
we use the equivalence with the Euler-Poincar\'e approach obtained above to write the equations in a simpler form.

Since $L_{ m _0 }$ is $H$-invariant, by fixing a connection $ \mathcal{A} \in \Omega ^1 (G, \mathfrak{h}  )$, we get the Lagrange-Poincar\'e Lagrangian $ \mathcal{L} $, that we now compute. We have
\[
\mathcal{L} \left( m, \dot m, \left[ g, \mathcal{A} (g, \dot g)\right] _H \right)= L_{ m _0 }(g, \dot g)= \mathfrak{L}  ( g m _0 , (\dot g g ^{-1} )_M( g m _0 ))=\mathfrak{L} (m, \dot m).
\]
This means that $\mathcal{L} :  TM\times_M\widetilde{\mathfrak{h}} \rightarrow \mathbb{R}$ 
does not depend on the second variable, so 
$ \frac{\delta \mathcal{L} }{\delta \sigma }=0$, and $\mathcal{L}  = \mathfrak{L} $. Thus, in the general system \eqref{LP_equations}, the first equation disappears and the right hand side of the second vanishes. Therefore, the Lagrange-Poincar\'e equations in \eqref{LP_equations} 
reduce to the Euler-Lagrange equations for $\mathfrak{L}$ on $TM$. 

We now compute the Euler-Poincar\'e reduced Lagrangian. We have
\[
l(\xi, m)=l\left( \dot g  g^{-1} , \Phi _{g}( m _0) \right)= L_{ m _0 } (g, \dot g)= \mathfrak{L}  ( g m _0 , (\dot g g ^{-1} )_M( g m _0 ))=\mathfrak{L} (m,  \xi_M(m)).
\]
From the above results, we know that the Euler-Lagrange equations for 
$L_{ m _0 }$ are equivalent to the Euler-Poincar\'e equations
\[
\frac{d}{dt}  \frac{\delta l}{\delta \xi }+
\operatorname{ad}^*_\xi  \frac{\delta l}{\delta \xi }= 
\mathbf{J} \left(\frac{\delta l}{\delta m} \right), \quad  \dot m = \xi _M (m).
\]

From this discussion together with the facts recalled in 
\S\ref{sec:general_theory}, we obtain the following fundamental 
result, to be used in the rest of the paper.

\begin{theorem}
\label{EP_reduction_corollary}
The following statements are equivalent.
\begin{itemize}
\item[{\rm (i)}] The curve $m:[ t _0 , t _1 ] \rightarrow M$ is a solution of the Euler-Lagrange equations for $\mathfrak{L} : TM
\rightarrow \mathbb{R}$, i.e.,
\begin{equation*}
\frac{\partial  \mathfrak{L} }{\delta m}- \frac{d}{dt} 
\frac{\partial  \mathfrak{L} }{\delta \dot m}=0.
\end{equation*}
\item[{\rm (ii)}] The curve $m:[ t _0 , t _1 ] \rightarrow M$ is a solution of the Euler-Poincar\'e equations for 
$l:\mathfrak{g} \times M \rightarrow \mathbb{R}$, i.e.,
\begin{equation}\label{cov_EP} 
\frac{d}{dt} \frac{\delta l}{\delta \xi }+
\operatorname{ad}^*_\xi  \frac{\delta l}{\delta \xi }= 
\mathbf{J} \left(\frac{\delta l}{\delta m} \right), 
\quad  \dot m= \xi _M (m).
\end{equation}
\end{itemize}
\end{theorem}

\paragraph{Ginzburg-Landau theory of phase transitions.} 
We briefly review the major steps in Landau's theory of phase transitions.

Phenomenological Ginzburg-Landau theory, initially formulated
to describe the behavior of superconductivity and superfluidity of ${}^4$He near points of phase transition,
turned out to be also very convenient in the determination of of phase transitions of superfluid ${}^3$He. We recall here
briefly the main statements of Landau's second order
theory of phase transitions. For a detailed presentation,
see \cite[\S83, \S141--153, \S162]{LaLi1980}.

\noindent \textbf{1.)} At a phase transition point, the 
symmetry of the system spontaneously changes.

\noindent \textbf{2.)} The system is characterized by some macroscopic quantity, an \textit{order parameter}, e.g., 
the director field, the $Q$-tensor, or the wryness tensor
in liquid crystals. 

\noindent \textbf{3.)} Near the transition point, due to the
smallness of the parameter $\alpha_i(T-T_c)$, the 
free energy (i.e., the thermodynamic potential) admits
an expression of the following type
\[
F(p, T) = F_0(p) + \alpha_1(T-T_c) Q(\varphi^2) + 
\alpha_2(T-T_c) Q(\varphi^2)Q(\varphi^4) + 
\|\operatorname{grad} \varphi\|^2,
\]
where $Q(\varphi^2)$, $Q(\varphi^4)$ are invariant under
the symmetry group of a system of second and fourth order,
$T$ is the temperature, $T_c$ is the critical temperature 
at which the phase transition occurs, $p$ is the pressure, 
and $\varphi$ is an order parameter of the given physical
system which is chosen by the concrete physical situation
under study.

\noindent \textbf{4.)} The change of symmetry in the 
transition is determined only by the order parameter.

\noindent \textbf{5.)} It is possible to ignore fluctuations
of the oder parameter beyond $(T_c/q_c)^4$. The 
Levanyuk-Ginzburg criterion ensures the validity of the
expression of $F(p,T)$ given above, if the mean square fluctuation
of the parameter $\varphi$, averaged over the correlation volume,
is small compared with the characteristic value of $\left\langle \varphi \right\rangle$ (see \cite{LaLi1980}).

The advantage of the Ginzburg-Landau approach is based on 
the fact that, with relatively few basic assumptions, it
is possible to reduce the investigation (in many important
cases) of an infinite dimensional quantum particle system
to the study of a finite dimensional mechanical problem.
Superfluid ${}^3$He provides such an example. Of course,
this is a more complicated system than superfluid ${}^4$He
since there are more thermodynamic phases.

We describe now the concrete method implementing this
Ginzburg-Landau phase transition theory. One is given
a free energy, the sum of a potential $U(A)$ depending 
on some order parameter $A$, but not on its spatial 
derivatives, and a gradient term $F_{grad}(A, \nabla A)$ that 
depends on both the order parameter $A$ and its derivatives 
$\nabla A$.

\textbf{(A)} Expand the potential $U(A)$ up to fourth order
and replace it with this expression.

\textbf{(B)} Find the largest possible Lie group
that leaves this fourth degree polynomial $U(A)$ invariant
and determine the Lie group action (very often a representation) on the space of all order parameters $A$.

\textbf{(C)} Find the formal minima of the potential
$U(A)$, i.e., $\frac{\delta U}{\delta A}= 0$ and 
$\frac{\delta^2 U}{\delta A^2}\geq 0$ (positive Hessian).

\textbf{(D)} Take the Lie group orbit through each minimum
and consider it as a configuration space of a Lagrangian
system given by the free energy. Note that it is not 
necessary to add the potential $U(A)$ to the gradient term, 
since it is constant on each such orbit, by construction. 
The goal is the study of each Lagrangian system on such 
a Lie group orbit, because the Ginzburg-Landau equations 
turn out to be the Euler-Lagrange equations for $F_{grad}$.
Often, $F_{grad}$ determines a metric on the orbit.

The point is that, very often, the Lie group orbits of 
interest are finite dimensional, whereas the original
problem, whose total Lagrangian is the sum of the free
energy $F$ and the potential $U(A)$, is an infinite 
dimensional problem. Working on such orbits reduces hence
the given infinite dimensional problem to a finite dimensional one.

Sometimes, steps (B) and (C) are hard to carry out. In
practice one starts with a Lie group that is, on physical
grounds, a symmetry of the system and then determines, using
invariant theory, the most general polynomial of fourth
degree, invariant under this group. This polynomial is then taken 
as the potential $U(A)$. Then one classifies 
all orbits of this Lie group on the space of all order
parameters $A$ (or, at least, determines enough orbits) and finds, in this way, orbits of physical
interest that describe different thermodynamic phases of the system. This problem is solved using techniques developed in \cite{GoMo1978a, GoMo1978b}, \cite{BoMo1987}, \cite{Mo1993},
where thermodynamic phases are identified with orbits
containing a minimum of the potential $U(A)$ of the free energy.

It turns out that different types of 
textures for the system are given as solutions to the
Ginzburg-Landau equations for a given phase and that,
on each Lie group orbit, the Ginzburg-Landau equations 
are the Euler-Lagrange equations for $F_{grad}$.

We shall apply this method to the study of 
different phases in superfluid ${}^3$He 
(\cite{BoMo1987}) and rotating neutron stars
 (\cite{MoSa2011}). Using the same techniques one
can also study  
one-dimensional textures in liquid crystals and 
superfluids as well as  phase transitions between biaxial 
and uniaxial nematics; we leave these latter topics for a future
publication.

\subsection{One-dimensional textures in the A-phase of liquid 
Helium ${}^3$He}

The order parameter of superfluid ${}^3$He is given by complex 
$3 \times 3$ matrices $A \in \mathfrak{gl}(3,\mathbb{C})$\footnote{The matrices $A$ are related to 
$\Delta'_{\sigma \sigma'}$, the ``energetic gap'' of the triplet
pairing of interacting quasiparticles of ${}^3$He, and so this
gap can be expressed in terms of $A$. Thus $A$ can be regarded
as the order parameter of superfluid ${}^3$He; see 
\cite[\S5.2.1]{Mo1993}.}.

The free energy is given by
\[
\mathcal{F} (A, \nabla A)= F_{grad}(A , \nabla A) +U(A),
\]
where 
\[
F_{grad}(A , \nabla A)=  \gamma_1 
\sum_{i,p,k}\left(\partial_k\bar A_{pi}\right) 
\left(\partial_k A_{pi}\right)  + \gamma_2 \sum_{i,p,k}
\left(\partial_k \bar A_{pi} \right)
\left(\partial_i A_{pk} \right) +\gamma_3 \sum_{i,p,k}
\left(\partial_k \bar A_{pk} \right)
\left(\partial_i A_{pi} \right),
\]
$ \gamma _1 , \gamma _2 , \gamma _3 >0$ are constants, and $U(A)$ is in the Ginzburg-Landau form, namely,
\begin{align*} 
U(A)&= \alpha \operatorname{Tr}(AA^*) + 
\beta_1 |\operatorname{Tr}(AA^T)|^2  + \beta_2\left[\operatorname{Tr}(AA^*)\right]^2 + 
\beta_3\operatorname{Tr} \left[(A^*A)\overline{(A^*A)}\right]
\\
& \qquad + \beta_4\operatorname{Tr}\left[(AA^*) ^2\right]
+ \beta_5\operatorname{Tr}\left[(AA^*) \overline{(AA^*)}
\right],
\end{align*}
for $ \alpha , \beta _1 ,..., \beta _5 \in \mathbb{R}$. Note that these expressions are real valued.

In one dimension, we compute
\begin{equation}\label{Fgrad_1D}
\begin{aligned} 
F_{grad}(A , \partial _z A)&= \gamma _1 \partial _z \bar A_{pi} \partial _z A_{pi}+ \gamma _2 \partial _z \bar A_{p 3} \partial _z A_{p3}+ \gamma _3 \partial _z\bar  A_{p3} \partial _z A_{p3}\\
&=\operatorname{Re} \operatorname{Tr}( \Gamma \partial _z A ^\ast \partial _z A)= \left\langle \! \left\langle \partial _z A, \partial _z A \right\rangle \! \right\rangle,
\end{aligned}
\end{equation}  
where $ \Gamma = \operatorname{diag}( \gamma _1 , \gamma _1 , \gamma _1 + \gamma _2 + \gamma _3 )$ and we defined the inner product on $ \mathfrak{gl}(3, \mathbb{C}  )$ by 
\begin{equation}\label{inner_product} 
\left\langle \! \left\langle A, B \right\rangle \! \right\rangle := \operatorname{Re} \operatorname{Tr}(\Gamma A^\ast B ),\quad \Gamma := \operatorname{diag}( \gamma _1 , \gamma _1 , \gamma _1 + \gamma _2 + \gamma _3 ).
\end{equation} 
The following identities are useful in the computations:
\[
\left\langle \! \left\langle A, B \right\rangle \! \right\rangle = 
\left\langle \! \left\langle B, A \right\rangle \! \right\rangle, \qquad
\left\langle \! \left\langle u A, B \right\rangle \! \right\rangle = 
\left\langle \! \left\langle A, \overline{u} B 
\right\rangle \! \right\rangle,
\]
for any $A, B \in \mathfrak{gl}(3,\mathbb{C})$ and $u \in \mathbb{C}$.
In addition $\left\langle \! \left\langle\,, \right\rangle \! \right\rangle$ is $\mathbb{R}$-bilinear.

\paragraph{Group representation, orbits, and thermodynamic phases.}
The potential function $U(A)$ is invariant under the left representation of the compact Lie group $G=U(1) \times SO(3)_L \times SO(3)_R$ on 
$\mathfrak{gl}(3,\mathbb{C})$ given by
\begin{equation}
\label{general_action_so3}
\left(e^{{\rm i}\varphi}, R_1, R_2\right) \cdot A : = 
e^{{\rm i}\varphi}R_1AR_2^{-1},
\end{equation}
where $A \in \mathfrak{gl}(3,\mathbb{C})$ and $\left(e^{{\rm i}\varphi}, R_1, R_2\right) \in G$. As the formula above shows,
the indices $L$ and $R$ on the two groups $SO(3)$ indicate the
side of the multiplication on the matrix $A$.

Note that the term $F_{grad}$ is not $G$-invariant. However, to determine the thermodynamic phases, it suffices to study $U(A)$. The phases correspond to different orbits. A partial classification of the orbits is given 
in \cite{BoMo1987}. Below we shall consider only some of these orbits that are physically relevant for the phases of superfluid 
${}^3$He (see \cite[\S5.2]{Mo1993}).

The A-phase of superfluid ${}^3$He has two regimes depending on whether
$L \ll L_{dip}$ or $L \gg L_{dip}$, where $L$ and $L_{dip}$ are the
characteristic and dipole length, respectively. The first regime
corresponds to minimal degeneracy and the dipole interaction can be
neglected. The order parameter matrix $A \in \mathfrak{gl}(3,\mathbb{C})$ is representable as an element of the 
$U(1) \times SO(3)_L \times SO(3) _R$-orbit through the point $A_0$ given in \eqref{A_0_A_Phase}, i.e., 
$A= e^{{\rm i} \varphi }R _1 A _0 R _2 ^{-1} $. In the second regime, 
the energy of the dipole interaction should be taken into account. As 
a consequence, the order parameter matrix $A \in \mathfrak{gl}(3,\mathbb{C})$ is representable as an element of the 
$SO(3)$-orbit through the same matrix $A_0$ under the
different action $A= R A _0 R ^{-1}$. For details, see 
\cite[\S 5.2.3]{Mo1993}.

\subsubsection{The $A$-phase -- first regime}

We consider the orbit $M$ of $U(1) \times SO(3)_L \times SO(3) _R$ through the point
\begin{equation}\label{A_0_A_Phase} 
A _0 = \left(
\begin{array}{ccc}
0&0&0\\
0&0&0\\
1&{\rm i}&0
\end{array} 
\right)\in \mathfrak{gl}(3, \mathbb{C}).
\end{equation} 
We note that $e^{ {\rm i}\varphi}A_0 = \rho(\varphi)A_0 \rho(-\varphi)$,
where $\rho(\varphi) : = \exp(\varphi\widehat{\mathbf{e}}_3)$.

\begin{proposition} 
{\rm (i)} The isotropy subgroup of $A_0 $ is 
\[
H=\{(e^{{\rm i}\varphi}, \rho (\alpha) J_+, \rho(\varphi)\tilde J_+),(e^{{\rm i}\varphi},\rho(\alpha) J_-, \rho(\varphi )\tilde J_-)\} 
\subset G=U(1) \times SO(3)_L \times SO(3)_R,
\]
where $\rho(\alpha)= \operatorname{exp}(\alpha\widehat{\mathbf{e}}_3)$ and
\begin{equation}
\label{J_tilde_J}
J_\pm=  \left(
\begin{array}{ccc}
1&0&0\\
0&\pm 1&0\\
0&0&\pm 1
\end{array} 
\right),\quad  \tilde J_\pm=\left(
\begin{array}{ccc}
\pm 1&0&0\\
0&\pm 1&0\\
0&0&1
\end{array} 
\right).
\end{equation}
{\rm (ii)} We have the diffeomorphism
\begin{equation}\label{diffeo1} 
G/H \ni [e^{{\rm i}\varphi }, R _1, R _2 ]_H\longmapsto 
[R_2 \rho (-\varphi), R_1\mathbf{e}_3]_{\mathbb{Z}_2 } \in 
(SO(3) \times S^2 )/ \mathbb{Z}_2,
\end{equation} 
where $ \mathbb{Z}  _2 =\{ \pm 1\}$ acts on $(A , \mathbf{x} ) \in SO(3) \times S ^2 $ as $(-1) \cdot ( A, \mathbf{x} )= (A \tilde{J}_- , - \mathbf{x} )$.\\
{\rm (iii)} We have the diffeomorphism
\begin{equation}\label{diffeo2} 
(SO(3) \times S ^2 )/ \mathbb{Z}  _2\ni [ A, \mathbf{x} ]_{ \mathbb{Z}  _2 } \longmapsto \mathbf{x} \otimes ( A _1 + {\rm i} A _2 ) \in \operatorname{Orb}(A _0 ),
\end{equation} 
where $A _i $ denotes the $i^{th}$ column of the matrix $A$.
\end{proposition} 
\textbf{Proof.} {\rm (i)} Writing $A_0=\Re(A_0)+ {\rm i}\Im(A_0)$, where $\Re(A_0)$ and $\Im(A_0)$ are real and imaginary parts of $A _0 $, the equality $e^{{\rm i}\varphi}R_1A_0 R_2^{-1}=A_0$ is equivalent to the two equations
\begin{align*} 
&(\cos\varphi) R_1 \Re(A_0) R_2^{-1} - (\sin\theta) R_1 \Im(A_0) R_2^{-1} = \Re(A_0) \\
&(\sin\varphi) R_1 \Re(A_0) R_2^{-1} +
(\cos\theta) R_1 \Im(A_0) R_2^{-1} = \Im(A_0).
\end{align*}
The proof then follows from a direct computation that is done by writing the matrices $R _i $ in terms of their rows.\\
{\rm (ii)} Let us first show that the map is well-defined. Given $(e^{{\rm i}\varphi }, R _1, R _2 ) \in G$, any element in the equivalence class $[(e^{{\rm i}\varphi }, R _1, R _2 )]_H$ has the form 
\begin{equation}\label{another_representative} 
(e^{{\rm i}\varphi}, R_1, R_2 )(e^{{\rm i}\psi}, \rho(\alpha) J_\pm, 
\rho(\psi)\tilde{J}_\pm)= (e^{{\rm i} (\varphi + \psi)}, 
R_1 \rho (\alpha)J_\pm, R_2 \rho (\psi)\tilde{J}_\pm),
\end{equation} 
where $(e^{{\rm i}\psi}, \rho(\alpha) J_\pm, \rho(\psi)\tilde{J}_\pm)
\in H$. Applying formula \eqref{diffeo1} to the 
expression \eqref{another_representative}, yields
\begin{align*}
[R_2 \rho ( \psi )\tilde{J}_\pm \rho (- \varphi ) \rho (- \psi ), R _1 \rho ( \alpha )J_\pm \mathbf{e} _3 ]_{ \mathbb{Z}  _2 }= [ R _2  \rho (- \varphi ) \tilde{J}_\pm, \pm R _1 \mathbf{e} _3 ]_{ \mathbb{Z}  _2 }= [ R _2  \rho (- \varphi ) ,R _1 \mathbf{e} _3 ]_{ \mathbb{Z}  _2 },
\end{align*} 
where we used the properties
\begin{equation}\label{two_properties} 
\rho ( \alpha )\tilde{J}_\pm= \tilde {J}_\pm \rho ( \alpha ) \quad\text{and}\quad \rho (\alpha) J_\pm=J_\pm \rho (\pm \alpha ).
\end{equation} 
The map is clearly surjective. To show the injectivity, we take 
$(e^{{\rm i}\varphi}, R_1, R_2), (e^{{\rm i}\varphi '}, R_1', R'_2 )  
\in G$ such that $[R_2 \rho (-\varphi), R_1\mathbf{e}_3]_{\mathbb{Z}_2}
=[R_2' \rho (-\varphi'), R_1'\mathbf{e}_3]_{\mathbb{Z}_2 }$. 
We thus have the equalities $R_2 \rho (-\varphi )=R_2'\rho (-\varphi ')\tilde{J}_\pm$ 
and $R_1 \mathbf{e}_3 =\pm R_1' \mathbf{e}_3 $. From the first 
equality, there exists $\rho (\psi) \in U(1)$ such that 
$R_2=R_2' \rho (\psi)\tilde{J}_\pm$, by using 
\eqref{two_properties}. Rewriting the second equality as 
$R_1 \mathbf{e}_3 = R_1 ' J_\pm \mathbf{e}_3 $, we obtain the 
existence of $\rho (\alpha) \in U(1)$ such that $R_1=
R_1' \rho (\alpha)J_\pm$. This proves that $[(e^{{\rm i}\varphi }, 
R_1, R _2 )]_H=[(e^{{\rm i}\varphi '}, R_1', R'_2 )]_H$.\\
{\rm (iii)} Given $[A, \mathbf{x}] _{ \mathbb{Z} _2 }\in (SO(3) \times S ^2 )/ \mathbb{Z}  _2$, let $(e^{{\rm i}\varphi }, R _1, R _2 ) \in G$ be such that $[R_2 \rho (-\varphi), R_1\mathbf{e}_3]_{\mathbb{Z}_2 }= [A, \mathbf{x} ]_{ \mathbb{Z}  _2 }$. A possible choice is $(e^{{\rm i}\varphi }, R _1, R _2 )=(1, R _1 , A)$, where $R _1 \in SO(3)$ is such that $ R _1 \mathbf{e} _3 = \mathbf{x} $. With this choice, an easy computation shows that
\[
(1, R _1 , A)\cdot A _0 = \mathbf{x} \otimes ( A _1 + {\rm i} A _2 ).\qquad\blacksquare
\]

\begin{remark}
\label{g_tilde_rem}
{\rm
We observe that the subgroup $\tilde G: = SO(3)_L \times SO(3)_R \subset G$ acts transitively on the orbit $ \operatorname{Orb}(A_0)$, see \eqref{diffeo1}, \eqref{diffeo2}. The isotropy subgroup of $A_0$ is $\tilde G_{A _0 }=H \cap \tilde G=\{1, \rho ( \alpha )J, \tilde J\}$, which is isomorphic to $O(2)$. Therefore the orbit can be equally well described as the homogeneous space 
$(SO(3)_L \times SO(3)_R)/\tilde G_{A _0}$.\quad $\blacklozenge$}
\end{remark}

\paragraph{Lagrangian formulation.}
We now apply Theorem \ref{EP_reduction_corollary} with this description of the orbit, so the Lie algebra is $ \mathfrak{so}(3)_L \times \mathfrak{so}(3)_R$.
On this orbit $M$, we consider the Lagrangian density given by the gradient part only, i.e.,
\begin{equation}\label{Lagrangian_A_phase}
\mathcal{L}(A, \partial _z  A) = F_{grad}(A, \partial _z  A)= \left\langle (\partial _z A) \Gamma , \partial _z A \right\rangle =\left\langle \! \left\langle \partial _z A, \partial _z A \right\rangle \! \right\rangle,
\end{equation} 
where
\begin{equation}\label{pairing}
\left\langle A, B \right\rangle := 
\operatorname{Re} \operatorname{Tr}(A ^\ast B).
\end{equation}
Note that 
\begin{equation}
\label{partial_derivative}
\frac{\delta \mathcal{L}}{\delta \partial_zA} = 2 \partial_z A \Gamma.
\end{equation}
The texture equations are given by the Euler-Lagrange equations
for $\mathcal{L}$ on the orbit $M$.

The reduced velocity $\xi =\partial _z g g^{-1}$ of the general theory (see 
\eqref{reduced_section_formula_for_g}) is given 
here by $ \xi = (\mathbf{v} ,\mathbf{w} ) : \mathbb{R} \rightarrow  \mathfrak{so}(3)_L \times \mathfrak{so}(3) _R$, where $\mathbf{v}$ and $\mathbf{w} $ are the chiral velocities $\mathbf{v} = (\partial_z R_1) R_1^{-1}$ and 
$\mathbf{w} = R_2^{-1} (\partial_z R _2)$, see \cite[formula (5.133)]{Mo1993}. The second formula in \eqref{EP_equations} is given here by
\[
\partial _z A= \widehat{ \mathbf{v} }A+ A\widehat{ \mathbf{w} }.
\] 
Using this expression and formula \eqref{Fgrad_1D}, the Euler-Poincar\'e Lagrangian
\[
l=l(\xi, m):  \mathfrak{so}(3)_L\times \mathfrak{so}(3)_R  \times M \rightarrow 
\mathbb{R}
\]
of the general theory given in  \eqref{EP_Lagrangian}, is computed in this case to be
\begin{align*} 
l(\mathbf{v} ,\mathbf{w} ,A) &= \operatorname{Re} \operatorname{Tr}( \Gamma \partial _z A ^\ast \partial _z A)= \left\langle \! \left\langle \widehat{ \mathbf{v} }A+ A\widehat{ \mathbf{w} }, \widehat{ \mathbf{v} }A+ A\widehat{ \mathbf{w} } \right\rangle \! \right\rangle ,
\end{align*} 
(see \eqref{inner_product}). Defining
\[
I_{ab}(A)= \left\langle\!\left\langle A \widehat{\mathbf{e} }_a, A \widehat{\mathbf{e} } _b \right\rangle \! \right\rangle , \quad \chi _{ab}(A)= \left\langle \! \left\langle \widehat{\mathbf{e} } _a A,\widehat{\mathbf{e} } _b A\right\rangle \! \right\rangle, \quad\text{and}\quad \Sigma _{ab}(A)= \left\langle\!\left\langle  \widehat{\mathbf{e} }_a A, A \widehat{\mathbf{e} } _b \right\rangle \! \right\rangle,
\]
the formula for the Lagrangian above becomes
\begin{align*} 
l(\mathbf{v} ,\mathbf{w} ,A)& = \sum_{a,b=1}^3\left( I_{ab}(A) w _a w _b + \chi _{ab}(A) v _a v _b+2 \Sigma _{ab}(A)v _a w _b \right)\\
&= \mathbf{w} ^\mathsf{T}\mathbf{I} (A)\mathbf{w} + \mathbf{v} ^\mathsf{T} \boldsymbol{\chi} (A) \mathbf{v} + 2 \mathbf{v} ^\mathsf{T} \boldsymbol{\Sigma} (A) \mathbf{w}.
\end{align*} 
Thus,
the Euler-Poincar\'e equations \eqref{EP_equations} read
\begin{equation*}
\partial _z \frac{\delta l}{\delta \mathbf{v} } + 
\operatorname{ad}^*_ \mathbf{v}  \frac{\delta l}{\delta \mathbf{v} }= 
\mathbf{J}_1   \left( \frac{\delta l}{\delta A} \right), \quad \partial _z \frac{\delta l}{\delta \mathbf{w} } - 
\operatorname{ad}^*_ \mathbf{w}  \frac{\delta l}{\delta \mathbf{w} }= 
\mathbf{J} _2  \left( \frac{\delta l}{\delta A} \right) ,\quad \partial _z A= \widehat{ \mathbf{v} }A+ A\widehat{ \mathbf{w} },
\end{equation*}
where $ \mathbf{J}_1 :T^*M \rightarrow \mathfrak{so}^*(3)$ is 
the momentum map of the left action and $\mathbf{J}_2 :T^*M 
\rightarrow \mathfrak{so}^*(3)$ is the momentum map of the right action of $SO(3)$ on the orbit $M$, respectively.

Using the duality pairing $ \left\langle A, B \right\rangle = \operatorname{Re} \operatorname{Tr}(A ^\ast B)$ on $ \mathfrak{gl}(3, \mathbb{C}  )$, we get the Euler-Poincar\'e equations
\[
\frac{d}{dz} \frac{\delta l}{\delta \mathbf{v} }+  \frac{\delta l}{\delta \mathbf{v} } \times \mathbf{v} = 2\overrightarrow{\operatorname{Re}\left(  \frac{\delta l}{\delta A}  A ^\ast \right)^{\phantom{4}}\!\! }, \qquad \frac{d}{dz} \frac{\delta l}{\delta \mathbf{w} }- \frac{\delta l}{\delta \mathbf{w} } \times \mathbf{w} =2\overrightarrow{\operatorname{Re}\left(   A ^\ast \frac{\delta l}{\delta A} \right)^{\phantom{4}}\!\! },
\]
where $\overrightarrow{A^{\phantom{4}}\!\!}\in \mathbb{R}  ^3 $ is defined by $\widehat{\overrightarrow{A^{\phantom{4}}\!\!}}:= A^{skew}:= \frac{1}{2} (A-A^\mathsf{T})$, and where we have
\[
\frac{\delta l}{\delta \mathbf{v} }=2 \boldsymbol{\chi} \mathbf{v}+2 \boldsymbol{\Sigma} \mathbf{w} ,  \quad \frac{\delta l}{\delta \mathbf{w} }=2 \mathbf{I} \mathbf{w}+2 \boldsymbol{\Sigma} ^\mathsf{T}\mathbf{w}, \quad \frac{\delta l}{\delta A}= -2 \left(  A\widehat{ \mathbf{w} } \Gamma \widehat{\mathbf{w} }+ \widehat{ \mathbf{v} }\widehat{ \mathbf{v} }A \Gamma+ \widehat{ \mathbf{v} }A \widehat{ \mathbf{w} }\Gamma + \widehat{ \mathbf{v} }A\Gamma \widehat{ \mathbf{w} } \right).
\]

\paragraph{Hamiltonian formulation.} As expected from the general theory, the Euler-Poincar\'e Lagrangian $l( \mathbf{v} , \mathbf{w} , A)$ is degenerate, since for all $A \in M$, the quadratic form $( \mathbf{v} , \mathbf{w} ) \mapsto\left\langle \! \left\langle \widehat{ \mathbf{v} }A+ A\widehat{ \mathbf{w} }, \widehat{ \mathbf{v} }A+ A\widehat{ \mathbf{w} } \right\rangle \! \right\rangle$ has a one dimensional kernel given by the isotropy Lie algebra $ \mathfrak{g}  _A=\{ ( \mathbf{v}, \mathbf{w} ) \in \mathfrak{g}  \mid\widehat{ \mathbf{v} }A+ A\widehat{ \mathbf{w} }=0\}$.

Since the Lagrangian \eqref{Lagrangian_A_phase} is nondegenerate, we consider the associated Hamiltonian on $T^*M$, given by
\begin{equation}\label{Hamiltonian_A_Phase} 
\mathcal{H} ( \alpha _A )= \frac{1}{4} \operatorname{Re} \operatorname{Tr}( \Gamma ^{-1} \alpha _A ^\ast \alpha _A )  = \frac{1}{4} \left\langle \alpha _A \Gamma ^{-1} , \alpha _A \right\rangle.
\end{equation} 

Now, we apply Theorem \ref{symplectic_leaves} in this particular case. The element $a _0 $ is given by $A_0$ in \eqref{A_0_A_Phase}. The groups are $G= SO(3)_L \times SO(3)_R$, $H=\widetilde{G}_{ A _0 }=\{ \rho ( \alpha ) J, \tilde{J}\}$. Given $ \mu =( \mathbf{m} , \mathbf{n} )\in \mathfrak{so}(3) ^\ast  \times \mathfrak{so}(3) ^\ast = \mathbb{R}  ^3 \times \mathbb{R}  ^3 $, since $ \mathfrak{g}  _{A _0 }=\{( \lambda \mathbf{e} _3 , \mathbf{0}) \mid \lambda \in \mathbb{R}  \}$, we have $ \mu _{ a _0 }=( m _3 \mathbf{e} _3 , \mathbf{0} )$. We now compute the $G_{A _0 }$-coadjoint orbit $ \mathcal{O} _{ \mu _{ a _0 }}$. We have the formulas
\begin{align*} 
\operatorname{Ad}_{(\rho ( \alpha )J, \tilde J)}(\widehat{\mathbf{v} } , \widehat{\mathbf{w} })&= ( \rho ( \alpha ) J \widehat{ \mathbf{v}} J \rho (- \alpha ), \tilde J \widehat{\mathbf{w} } \tilde J)\\
\operatorname{Ad}^\ast_{(\rho ( \alpha )J, \tilde J)}(m _3 \mathbf{e} _3 , \mathbf{0} )&= ((-1)^{|J|} m _3 \mathbf{e} _3 , \mathbf{0} ),
\end{align*} 
where $|J|= 0$ if $J=I_3$ and $|J|=1$ otherwise. Thus, $ \mathcal{O} _{ (m_3 \mathbf{e} _3 , \mathbf{0} )}= \{(\pm m_3 \mathbf{e} _3 , \mathbf{0} )\}$ and hence the fibers of the associated fiber bundle $\widetilde{\mathcal{O}}_{ (m_3 \mathbf{e} _3 , \mathbf{0} )}\rightarrow M$  are two points sets.
In this special situation, the symplectic structure on $T^*M \times _M \widetilde{\mathcal{O}}_{ (m_3 \mathbf{e} _3 , \mathbf{0} )}$ is given by the canonical symplectic form on $T^*M$ since the Lie algebra  $ \mathfrak{g}  _{A _0 }$  is one-dimensional and the fiber is discrete, see \cite[Theorem 2.3.12]{MaMiOrPeRa2007}. We conclude that the coadjoint
orbit $\mathcal{O}_{(\mathbf{m}, \mathbf{n}, A_0)}$ has two connected
components each one symplectically diffeomorphic to $T^*M$ for any
$\mathbf{m}, \mathbf{n} \in\mathbb{R}^3$. In particular, the dimension
of the coadjoint orbit $\mathcal{O}_{(\mathbf{m}, \mathbf{n}, A_0)}$
is ten.

Now, we extend the Hamiltonian \eqref{Hamiltonian_A_Phase} to the symplectic manifold $T^*M \times _M \widetilde{\mathcal{O}}_{ (m_3 \mathbf{e} _3 , \mathbf{0} )}$. Hamilton's equations are unchanged. Using the symplectic diffeomorphism of Theorem \ref{symplectic_leaves} we get a Hamiltonian function on the coadjoint orbit $ \mathcal{O} _{( \mathbf{m} , \mathbf{n} , A _0)}$ of the semidirect product $(SO(3)_L\times SO(3)_R) \,\circledS\, \mathfrak{gl}(3, \mathbb{C}  )$. It is a symplectic leaf of the Lie-Poisson manifold $[(\mathfrak{so}(3)_L \times \mathfrak{so}(3)  _R) \,\circledS\, \mathfrak{gl}(3, \mathbb{C}  )] ^\ast $ and hence of its Poisson submanifold $(\mathfrak{so}(3) _L\times  \mathfrak{so}(3)  _R)^\ast  \times M$, endowed with the Lie-Poisson bracket
\begin{equation}\label{Lie_Poisson_bracket}
\begin{aligned} 
\{f,h\}( \mathbf{m} , \mathbf{n} , A)&= \mathbf{m} \cdot \frac{\delta f}{\delta \mathbf{m} }\times \frac{\delta h}{\delta \mathbf{m} }- \mathbf{n} \cdot \frac{\delta f}{\delta \mathbf{n} }\times \frac{\delta h}{\delta \mathbf{n} }\\
&+ \left\langle \frac{\delta f}{\delta A}, \widehat{ \frac{\delta  h}{\delta  \mathbf{m} }}A +A\widehat{ \frac{\delta  h}{\delta  \mathbf{n} }}\right\rangle-\left\langle \frac{\delta h}{\delta A}, \widehat{ \frac{\delta  f}{\delta  \mathbf{m} }}A +A\widehat{ \frac{\delta  f}{\delta  \mathbf{n} }}\right\rangle.
\end{aligned}
\end{equation} 

A direct computation shows that the kernel of the Poisson tensor 
is one dimensional at all points $(\mathbf{m},\mathbf{n},A_0)$. 
This means that the dimension of the symplectic leaves through 
$(\mathbf{m}, \mathbf{n}, A_0)$ is ten. We have recovered the 
previous result stating that the dimension of the coadjoint 
orbit $\mathcal{O}_{(\mathbf{m}, \mathbf{n}, A_0)}$ is ten. 
We note that the function $C(\mathbf{m}, \mathbf{n} , A)=
\frac{1}{2}  \operatorname{Re} \operatorname{Tr}(A^\ast A)$ 
is a Casimir function of this bracket. Indeed, since 
$\frac{\delta C}{\delta A}= A$, a direct computation that 
involves only the third term in the expression above shows 
that $\{C,f\}=0$ for all functions $f$.

\begin{lemma} The Riemannian metric on $M$ induced by the 
$\operatorname{Ad}$-invariant inner product 
$\gamma (( \mathbf{a}  , \mathbf{b} ), (\mathbf{v} , \mathbf{w} ))=\mathbf{a} \cdot \mathbf{v} + \mathbf{b} \cdot \mathbf{w} $ on $\mathfrak{so}(3)_L \times \mathfrak{so}(3)_R$ 
{\rm (}$M$ is viewed here as the orbit
$(SO(3)_L \times SO(3)_R)/\tilde G_{A _0}$ as in Remark \ref{g_tilde_rem}{\rm )} coincides with the metric induced by the inner product \eqref{pairing} {\rm(}here, $M \subset  
\mathfrak{gl}(3, \mathbb{C})${\rm)}, that is,
\[
\gamma _M( \widehat{ \mathbf{a} }A+A  \widehat{ \mathbf{b} },  \widehat{ \mathbf{v} }A+A  \widehat{ \mathbf{w} })=  \operatorname{Re} \operatorname{Tr}\left( (\widehat{ \mathbf{a} }A+A  \widehat{ \mathbf{b} }) ^\ast ( \widehat{ \mathbf{v} }A+A  \widehat{ \mathbf{w} })\right).
\]
\end{lemma} 
\textbf{Proof.} We need to verify identity \eqref{induced_metric}. 
It is readily checked that at $A_0$, we have $ \operatorname{Re} \operatorname{Tr} \left(    (\widehat{ \mathbf{a} }A_0 +A _0  \widehat{ \mathbf{b} }) ^\ast ( \widehat{ \mathbf{v} }A_0 +A_0   \widehat{ \mathbf{w} })\right)= a _1 v _1 + a _2 v _2 + \mathbf{b} \cdot \mathbf{w}= \mathbb{P}  _{A _0 }^\perp( \mathbf{a} , \mathbf{b} ) \cdot \mathbb{P}  _{ A _0 }^\perp( \mathbf{v} , \mathbf{w} )$, where $ \mathbb{P}_{A _0 }^\perp( \mathbf{a}, \mathbf{b} )= ((a _1 , a _2 , 0), \mathbf{b} )$.  Since $( \mathbf{a} , \mathbf{b} )_M(A)= \widehat{\mathbf{a}} A+A\widehat{\mathbf{b}}$, inserting the expression $A= R _1 A _0 R _2 ^{-1} $, we get
\begin{align*}
\operatorname{Re} \operatorname{Tr}\left( (\widehat{ \mathbf{a} }A+A  \widehat{ \mathbf{b} }) ^\ast ( \widehat{ \mathbf{v} }A+A  \widehat{ \mathbf{w} })\right)&=( R _1 ^{-1} \mathbf{a} )_1( R _1 ^{-1} \mathbf{v} )_1+( R _1 ^{-1} \mathbf{a} )_2( R _1 ^{-1} \mathbf{a} )_2+ R _2 ^{-1} \mathbf{b} \cdot R _2 ^{-1} \mathbf{w}\\
&=  \mathbb{P}_{A _0 }^\perp(R _1 ^{-1} \mathbf{a} , R _2 ^{-1} \mathbf{b} ) \cdot  \mathbb{P}_{A _0 }^\perp(R_1 ^{-1} \mathbf{b} , R _2 ^{-1} \mathbf{w}  ),
\end{align*} 
which proves the formula. $\quad\blacksquare$

\medskip

It follows that formula \eqref{Hor_Lift} can be applied. Therefore, we get
\[
\operatorname{Hor}^*_{( R _1 , R _2 )}(  \widehat{\mathbf{m} }R _1,R _2 \widehat{ \mathbf{n} }  )= \widehat{\mathbf{m} }A+ A  \widehat{\mathbf{n} }\in T^*_A M.
\]
Fixing $ \mu =(\mathbf{m} _0 , \mathbf{n} _0 )$ and applying Theorem \ref{symplectic_leaves} we get the Hamiltonian function on the coadjoint orbit $ \mathcal{O} _{( \mathbf{m}_0  , \mathbf{n}_0  , A _0 )}$ by pulling back the Hamiltonian $ \mathcal{H} $ in \eqref{Hamiltonian_A_Phase}. We obtain
\begin{equation}\label{Hamiltonian_coadjoint} 
h( \mathbf{m} , \mathbf{n} , A)= \mathcal{H} ( \widehat{\mathbf{m} }A+ A  \widehat{\mathbf{n} })= \frac{1}{4} \left\langle( \widehat{\mathbf{m} }A+ A  \widehat{\mathbf{n} } )\Gamma ^{-1} , \widehat{\mathbf{m} }A+ A  \widehat{\mathbf{n} } \right\rangle , 
\end{equation} 
where $( \mathbf{m} , \mathbf{n} , A) \in \mathcal{O} _{( \mathbf{m}_0  , \mathbf{n}_0  , A _0 )}$.

The general formula for the coadjoint action on a semidirect 
product \eqref{coadj_semidirect} (see, e.g., \cite{MaRaWe1984}) 
yields in this case
\[
\operatorname{Ad}^*_{( R _1 , R _2 , V)^{-1} }( \mathbf{m} , \mathbf{n} , A)= \left( R _1 \mathbf{m} +2 \overrightarrow{\operatorname{Re}(R _1 A R _2 ^{-1} V)},R _2 \mathbf{n} +2 \overrightarrow{\operatorname{Re}(VR _1 A R _2 ^{-1})}, R _1 A R _2 ^{-1}\right),
\]
where $( R _1, R _2 , V) \in (SO(3)_L \times SO(3)_R) \,\circledS\,  \mathfrak{gl}(3, \mathbb{C}  )$.

We will now consider subgroup actions of the coadjoint action that are symmetries of the Hamiltonian \eqref{Hamiltonian_coadjoint} and compute the associated momentum maps.

The first one is given by the $U(1)$-action $ \operatorname{Ad}^*_{(I_3, \rho ( \varphi ), 0)}( \mathbf{m} , \mathbf{n} , A)=   ( \mathbf{m} , \rho ( \varphi ) \mathbf{n} , A \rho (- \varphi ))$. This action is automatically Poisson and leaves the Hamiltonian \eqref{Hamiltonian_coadjoint} invariant because $ \rho ( \varphi ) \Gamma^{-1}  = \Gamma ^{-1} \rho ( \varphi )$. The infinitesimal generator of this action is $( \mathbf{m} , \mathbf{n} , A) \mapsto ( \mathbf{0}, \mathbf{e} _3 \times \mathbf{n} , -A \widehat{ \mathbf{e} }_3)$ and the momentum map is found to be $\mathbf{J}^{\rm orb}_ 3 ( \mathbf{m} , \mathbf{n} , A )=- \mathbf{e} _3 \cdot \mathbf{n}$. Therefore, $\{ \mathbf{J}^{\rm orb}_ 3, h\}=0$.

The second symmetry is given by the $SO(3)$-action $ \operatorname{Ad}^*_{(R, I_3, 0)}( \mathbf{m} , \mathbf{n} , A)=(R \mathbf{m} , \mathbf{n}, RA)$ whose infinitesimal generator associated to $\widehat{ \mathbf{v}  } \in \mathfrak{so}(3)  $ is $( \mathbf{m} , \mathbf{n} , A) \mapsto (\mathbf{v} \times \mathbf{m} , \mathbf{0}  ,\widehat{ \mathbf{v}  }A)$. This action leaves the Hamiltonian \eqref{Hamiltonian_coadjoint} invariant. The momentum map is $\mathbf{J}^{\rm spin} ( \mathbf{m} , \mathbf{n} , A )= \mathbf{m} $. Therefore $\{ \mathbf{J}^{\rm spin}_{ \mathbf{v} }, h\}=0$, for all $ \mathbf{v} \in \mathbb{R}  ^3 $. In particular, $\{\mathbf{J}^{\rm spin}_ 3, h\}=0$ and $\{\|\mathbf{J}^{\rm spin} \| ^2 , h\}=0$. In addition, formula \eqref{Lie_Poisson_bracket} implies that $\{\mathbf{J}^{\rm spin}_ 3, \mathbf{J}^{\rm orb}_ 3\}=0$ and $\{\|\mathbf{J}^{\rm spin} \| ^2,\mathbf{J}^{\rm orb}_ 3\}=0$.

To find the next conserved quantity is considerably more involved. We start with the Euler-Lagrange equations for the Lagrangian $ \mathcal{L} (A , \partial _z A)$ in \eqref{Lagrangian_A_phase} on $TM$. This Lagrangian is $U(1)$-invariant under the tangent lift of the action $ A \mapsto e^{i \varphi } A$. The infinitesimal generator associated to $ \theta \in \mathbb{R}  $ is $ \theta _M(A)= {\rm i} \theta A$. Using \eqref{partial_derivative}, we compute the associated momentum map as follows
\begin{equation}
\label{j_m_A_first}
j_m(A, \partial _z A) = \left\langle \frac{\delta  \mathcal{L} }{\delta  ( \partial _z A)},  {\rm i} A \right\rangle =2\left\langle \partial _z A \,\Gamma ,  {\rm i} A \right\rangle.
\end{equation}
Taking into account that $ \partial _z A= \widehat{\mathbf{v}}A+A\widehat{\mathbf{w}}$, for some $ \mathbf{v} , \mathbf{w} \in \mathbb{R}  ^3 $, this formula becomes
\[
j_m(A, \widehat{\mathbf{v}}A+A\widehat{\mathbf{w}})= -2 \operatorname{Re} \operatorname{Tr}(A \Gamma  A ^\ast \widehat{\mathbf{v}} {\rm i})-2 \operatorname{Re} \operatorname{Tr}( \Gamma \widehat{\mathbf{w}} A ^\ast  A{\rm i})=2 \left\langle \! \left\langle \widehat{\mathbf{v}} A +A \widehat{\mathbf{w}}, {\rm i}  A \right\rangle \! \right\rangle.
\]

Since $j _m $ is conserved on the solutions of the Euler-Lagrange equations associated to $ \mathcal{L} $, its pull-back to $ \mathcal{O} _{( \mathbf{m} _0 , \mathbf{n} _0 , A _0 )}$ commutes with the Hamiltonian $h$.

In order to see that $ j _m $ commutes with $\mathbf{J}^{\rm orb}_3$ 
and $\mathbf{J} ^{\rm spin}$, we will consider the induced $U(1)$ 
and $SO(3)$-actions on $TM$ and $T^*M$ and observe that they are 
the tangent and cotangent lift of commuting actions. Therefore, 
viewed as momentum maps on $T^*M$ and $TM$, via the change of 
variables $\mathcal{O}_{(\mathbf{m}_0 , \mathbf{n}_0 , A_0)} 
\rightarrow T^*M \rightarrow TM$ (see Theorem \ref{symplectic_leaves} 
and \eqref{partial_derivative}), these momentum maps commute.
Concerning $\mathbf{J}^{\rm orb}_3$, the $U(1)$-action induced on $TM$ is the tangent lift of the action $A \mapsto A\rho ( -\varphi )$.
For $ \mathbf{J} ^{\rm spin}$, the $SO(3)$-action induced on $TM$ is the tangent lift of the action $A\mapsto RA$, $R \in SO(3)$. They evidently commute with the action $ A \mapsto e^{i \varphi } A$ yielding $j _m $.

One can also check directly that the expressions of the momentum maps $j^{\rm orb}_3(A, \partial _z A)$ and $j^{\rm spin}(A, \partial _z A)$ on $TM$ associated to the tangent lifted actions of $A \mapsto A\rho ( -\varphi )$ and $A\mapsto RA$ are consistent with those of $\mathbf{J}^{\rm orb}_3( \mathbf{m} , \mathbf{n} , A)$ and $\mathbf{J}^{\rm spin}( \mathbf{m} , \mathbf{n} , A)$, respectively.

\begin{theorem} 
The five functions $h, j _m ,  \mathbf{J}^{\rm orb}_3, \mathbf{J} ^{\rm spin}_3, \|\mathbf{J} ^{\rm spin}\| ^2 $ form a completely integrable system on the ten dimensional coadjoint orbit $\mathcal{O} _{( \mathbf{m} _0 , \mathbf{n}_0 , A _0 )}$.
\end{theorem} 
\textbf{Proof}. The five functions commute in view of the discussion above. We need to show that their differentials are linearly independent except on a set of measure zero in $\mathcal{O} _{( \mathbf{m} _0 , \mathbf{n} _0 , A _0 )}$.
It turns out that showing their independence on $M$ is considerably simpler computationally.
The functional derivatives on $TM$ are
\begin{align*} 
&\frac{\delta j _m }{\delta A}= -2{\rm i} \partial _z A \, \Gamma , \quad \frac{\delta j _m }{\delta \partial _z A}=2 {\rm i}A \Gamma , \quad \frac{\delta j _3 ^{\rm orb}}{\delta A}= 2 \partial _z A \Gamma \widehat{ \mathbf{e} } _3, \quad \frac{\delta j^{\rm orb}_3 }{\delta \partial _z A}=-2 A \widehat{ \mathbf{e} } _3 \Gamma,\\
& \frac{\delta j _k ^{\rm spin}}{\delta A}= -2 \widehat{ \mathbf{e} } _k \partial _z A \, \Gamma , \quad   \frac{\delta j _k ^{\rm spin}}{\delta \partial _z A}= 2 \widehat{ \mathbf{e} } _k A  \Gamma, \quad \frac{\delta \mathcal{L} }{\delta A}=0, \quad \frac{\delta \mathcal{L} }{\delta \partial _z A}=  2 \partial _z A \,\Gamma,\\
& \frac{\delta \| \mathbf{J} ^{\rm spin}\| ^2 }{\delta A}=  - 4 j ^{\rm spin}_1 \widehat{ \mathbf{e} } _1 \partial _z A\, \Gamma - 4 j ^{\rm spin}_2 \widehat{ \mathbf{e} } _2 \partial _z A\, \Gamma, \quad  \frac{\delta \| \mathbf{J} ^{\rm spin}\| ^2 }{\delta \partial _z A}=  4 j ^{\rm spin}_1 \widehat{ \mathbf{e} } _1  A \Gamma +  4 j ^{\rm spin}_2 \widehat{ \mathbf{e} } _2  A \Gamma.
\end{align*} 
In order to show the independence, we have to show that the equations
\begin{align} 
&\alpha _1 \frac{\delta j _m }{\delta A}+ \alpha _2 \frac{\delta j _3 ^{\rm orb}}{\delta A}+ \alpha _3 \frac{\delta j _3 ^{\rm spin}}{\delta A}+ \alpha _4 \frac{\delta \| \mathbf{J} ^{\rm spin}\| ^2}{\delta A}+ \alpha _5 \frac{\delta \mathcal{L} }{\delta A}=0\label{first_equation} \\
&\alpha _1 \frac{\delta j _m }{\delta \partial _z A}+ \alpha _2 \frac{\delta j _3 ^{\rm orb}}{\delta \partial _z A}+ \alpha _3 \frac{\delta j _3 ^{\rm spin}}{\delta \partial _z A}+ \alpha _4 \frac{\delta \| \mathbf{J} ^{\rm spin}\| ^2}{\delta \partial _z A}+ \alpha _5 \frac{\delta \mathcal{L} }{\delta \partial _z A}=0\label{second_equation} 
\end{align}
imply $ \alpha _i =0$, for all $i=1,...,5$ and for all $A \in M$ except on a set of measure zero in $M$.

Writing $A= \mathbf{x}(\mathbf{A}_1 + {\rm i}\mathbf{A}_2)^\mathsf{T}$, 
where $\|\mathbf{x}\|=1$, $\|\mathbf{A}_i\|=1$, $\mathbf{A}_1 \cdot 
\mathbf{A}_2 =0$, and using the formula $\partial_z A= 
(\partial_z \mathbf{x})(\mathbf{A}_1 + {\rm i}\mathbf{A}_2)^\mathsf{T}
+ \mathbf{x}(\partial_z  \mathbf{A}_1 + {\rm i}\partial_z 
\mathbf{A}_2)^\mathsf{T}$, where $\partial_z \mathbf{x} \cdot \mathbf{x} 
=0$, $\partial_z \mathbf{A}_i \cdot \mathbf{A}_i =0$, $\partial_z 
\mathbf{A}_1 \cdot \mathbf{A}_2 + \mathbf{A}_1 \cdot \partial_z 
\mathbf{A}_2 =0$, and evaluating equation \eqref{first_equation} on 
the vector $(\mathbf{A}_1 + {\rm i}\mathbf{A}_2) \times 
(\partial_z \mathbf{A}_1 + {\rm i}\partial_z \mathbf{A}_2)$, we get 
$\alpha_2 (\partial_z \mathbf{A}_1 + {\rm i}\partial_z 
\mathbf{A}_2)^\mathsf{T}(\partial_z \mathbf{A}_1 + {\rm i}\partial_z \mathbf{A}_2) (\mathbf{e}_3 \cdot  (\mathbf{A}_1 + {\rm i}
\mathbf{A}_2))=0$. This implies $\alpha_2 =0$ except on a set of 
measure zero in $M$.

Using this and evaluating equation \eqref{second_equation} on 
$\mathbf{A}_1 \times \mathbf{A}_2 $, we get $\alpha_5 (\partial_z  
\mathbf{A}_1 + {\rm i}\partial_z\mathbf{A}_2 )^\mathsf{T}(\mathbf{A}_1 
\times \mathbf{A} _2 )=0$ which again implies $ \alpha _5 =0$ except on 
a set of measure zero in $M$. Then multiplying \eqref{second_equation} 
on the right by $A^\ast \mathbf{x}$, taking the dot product with 
$\mathbf{x}$, and using the formula $A A^\ast \mathbf{x} = 2\mathbf{x}$, 
we get $ \alpha_1 =0$. Multiplying the remaining equation on the left by 
$\mathbf{e}_3^\mathsf{T}$, we get $\alpha_4(j_1^{\rm spin}(A, 
\partial_z A)\mathbf{e}_2^\mathsf{T}-j_2^{\rm spin}(A, \partial_z A)
\mathbf{e}_1^\mathsf{T})\partial_z A =0$, which again, except on a set 
of measure zero in $M$, implies $\alpha_4 =0$. From this it follows 
that $\alpha_3 =0$ except on a set of measure zero in $M$. $\qquad\blacksquare$

\subsubsection{The $A$-phase -- second regime} 
\label{sec:A_phase_second}

In this situation, we consider the orbit $M=\{R A _0 R ^{-1} \mid R \in SO(3)\}$ of $SO(3)$ through $A _0 $ given by \eqref{A_0_A_Phase}.
A direct verification proves the following result.

\begin{proposition}
\label{prop:orbit_a_phase_second} 
${\rm (i)}$ The isotropy subgroup $SO(3)_{A_0}$ 
equals 
\[
SO(3)_{A_0}=\left\{I_3,\; 
\begin{pmatrix}
1&\;\;\,0&\;\;\;0\\
0&-1&\;\;\;0\\
0&\;\;\,0&-1
\end{pmatrix},\;
\begin{pmatrix}
-1&\;\;\,0&0\\
\;\;\,0&-1&0\\
\;\;\,0&\;\;\,0&1
\end{pmatrix},\;
\begin{pmatrix}
-1&\;\;\,0&\;\;\,0\\
\;\;\,0&\;\;\,1&\;\;\,0\\
\;\;\,0&\;\;\,0&-1
\end{pmatrix}
\right\} \cong \mathbb{Z}_2 \times \mathbb{Z}_2
\]
the group isomorphism being given by
\begin{align*}
I_3 &\longleftrightarrow (1,1),\qquad \qquad
\begin{pmatrix}
1&\;\;\,0&\;\;\;0\\
0&-1&\;\;\;0\\
0&\;\;\,0&-1
\end{pmatrix}
\longleftrightarrow (1,-1),\\ \\
\begin{pmatrix}
-1&\;\;\,0&0\\
\;\;\,0&-1&0\\
\;\;\,0&\;\;\,0&1
\end{pmatrix}
&\longleftrightarrow (-1,1),\qquad \quad
\begin{pmatrix}
-1&\;\;\,0&\;\;\,0\\
\;\;\,0&\;\;\,1&\;\;\,0\\
\;\;\,0&\;\;\,0&-1
\end{pmatrix}
\longleftrightarrow (-1,-1).
\end{align*}
${\rm (ii)}$ We have the diffeomorphism
\begin{equation*} 
SO(3)/SO(3)_{A_0} \ni [R]_{SO(3)_{A_0}}\longmapsto RA_0R^{-1} \in 
\operatorname{Orb}(A _0 ).
\end{equation*} 
\end{proposition} 

\paragraph{Lagrangian formulation.}
We apply Theorem \ref{EP_reduction_corollary} with this description of the orbit.  The reduced velocity (see 
\eqref{reduced_section_formula_for_g}) is given 
here by $\mathbf{w} = (\partial_z R )R ^{-1}$. The second formula in \eqref{EP_equations} becomes
\[
\partial _z A= \widehat{ \mathbf{w} }A- A\widehat{ \mathbf{w} }=[\widehat{ \mathbf{w} },A].
\] 
Using this expression and formula \eqref{Fgrad_1D}, the Euler-Poincar\'e Lagrangian \eqref{EP_Lagrangian} reads
\begin{align*} 
l(\mathbf{w} ,A) &= \operatorname{Re} \operatorname{Tr}( \Gamma \partial _z A ^\ast \partial _z A)= \left\langle \! \left\langle [A,\widehat{ \mathbf{w} }], [A,\widehat{ \mathbf{w} }] \right\rangle \! \right\rangle ,
\end{align*} 
(see \eqref{inner_product}). Defining
\begin{equation}
\label{J_A}
J_{ab}(A)= \left\langle\!\left\langle [A ,\widehat{\mathbf{e} }_a],[ A, \widehat{\mathbf{e} } _b] \right\rangle \! \right\rangle,
\end{equation}
the Lagrangian \eqref{Lagrangian_A_phase} reads
\[
l(\mathbf{w} ,A) = \sum_{a,b=1}^3 J_{ab}(A) w _a w _b = \mathbf{w} ^\mathsf{T}\mathbf{J} (A)\mathbf{w}.
\]
The Euler-Poincar\'e equations \eqref{EP_equations} are
\begin{equation*}
\partial _z \frac{\delta l}{\delta \mathbf{w} } + 
\operatorname{ad}^*_ \mathbf{w}  \frac{\delta l}{\delta \mathbf{w} }= -
\mathbf{J} \left( \frac{\delta l}{\delta A} \right),
\end{equation*}
where $ \mathbf{J} :T^*\operatorname{Orb}(A _0 ) \rightarrow \mathfrak{so}(3)^*$ is 
the momentum map of the right action of $SO(3)$ on $T^*\operatorname{Orb}(A _0 )$.
Using the duality pairing $ \left\langle A, B \right\rangle = \operatorname{Re} \operatorname{Tr}(A ^\ast B)$ on 
$\mathfrak{gl}(3,\mathbb{C})$, we get
\begin{equation}
\label{equ_A_phase_second}
\frac{d}{dz} \frac{\delta l}{\delta \mathbf{w} }+  \frac{\delta l}{\delta \mathbf{w} } \times \mathbf{w} = 2\operatorname{Re}\overrightarrow{\left[\frac{\delta l}{\delta A} ^\ast, A \right]^{\phantom{4}}\!\!},
\end{equation}
where we have
\[
\frac{\delta l}{\delta \mathbf{w} }=2 \mathbf{J}(A) \mathbf{\mathbf{w} }, \quad \frac{\delta l}{\delta A}= 2 [[\widehat{\mathbf{w} },A] \Gamma ,\widehat{\mathbf{w} }].
\]

\paragraph{Hamiltonian formulation.} Using the general formula 
\eqref{coadj_semidirect} for
the coadjoint action of the semidirect product $SO(3)\,\circledS\,
\mathfrak{gl}(3,\mathbb{C})$, it is it easily seen that the coadjoint
orbit through $(0, A_0)$ is $\mathfrak{so}(3)^\ast \times \operatorname{Orb}(A_0)$. This is consistent with the general
theory in Theorem \ref{symplectic_leaves}: since $\mathfrak{g}_{a_0}
= 0$, the dimension of the orbit is six.

Define $\mathbf{m}: = 
\frac{\delta l}{ \delta \mathbf{w}} = 2\mathbf{J}(A)\mathbf{w}$. Thus
the Hamiltonian associated to the Lagrangian $l$ has the
expression
\begin{equation}
\label{Hamiltonian_A_second}
h(\mathbf{m}, A) = \mathbf{m}^\mathsf{T}\mathbf{w} - l(\mathbf{w}, A) 
= \frac{1}{4} \mathbf{m}^\mathsf{T} \mathbf{J}(A)^{-1}\mathbf{m}.
\end{equation}
The non-degenerate Lie-Poisson bracket on the coadjoint orbit 
$\mathfrak{so}(3)^\ast \times \operatorname{Orb}(A_0)$ is
\begin{equation}\label{Lie_Poisson_bracket_A_second}
\begin{aligned} 
\{f,h\}(\mathbf{m} , A)&= \mathbf{m} \cdot \left(
\frac{\delta f}{\delta \mathbf{m}}\times 
\frac{\delta h}{\delta \mathbf{m}} \right)+ \left\langle \frac{\delta h}{\delta A}, \left[ A,\widehat{ \frac{\delta  f}{\delta  \mathbf{m} }}\right] \right\rangle-\left\langle \frac{\delta f}{\delta A}, \left[ A,\widehat{ \frac{\delta  h}{\delta  \mathbf{m} }}\right] \right\rangle
\end{aligned}
\end{equation} 
and hence the equations $\partial_z f = \{f, h\}$, for any $f$, are 
\begin{equation}
\label{LP_A_second_Ham_equ}
\partial_z\mathbf{m}+ \mathbf{m} \times 
\frac{\delta h}{\delta \mathbf{m}} = 
-2\operatorname{Re}\overrightarrow{\left[\frac{\delta h}{\delta A} ^\ast, A \right]^{\phantom{4}}\!\!}
\qquad 
\partial _z A=\left[\widehat{\frac{\delta h}{\delta \mathbf{m}}},A\right].
\end{equation}

We prove that the Hamiltonian system given by \eqref{Hamiltonian_A_second} relative to the Poisson bracket
\eqref{Lie_Poisson_bracket_A_second}  on the six dimensional
coadjoint obit $\mathfrak{so}(3)^\ast \times \operatorname{Orb}(A_0)$ is completely integrable. The three integrals of motion are the Hamiltonian
\eqref{Hamiltonian_A_second}, the momentum map $j_m$ given in
\eqref{j_m_A_first}, i.e., $j_m(\mathbf{m},A) = 
\left\langle \! \left\langle \left[\frac{1}{2} 
\widehat{\mathbf{J}(A)^{-1} \mathbf{m}}, A 
\right],{\rm i} A \right\rangle \! \right\rangle$ after transforming
to the variables $(\mathbf{m}, A)$, and $\mathbf{J}_3(\mathbf{m}, A): = \mathbf{e}_3 \cdot \mathbf{m}$. As in the discussion of the A-phase,
the previous regime, we note that $j_m$ is the momentum map associated 
to the circle action on configuration space given by $A \mapsto 
e^{{\rm i} \varphi}A$. Puling back $j_m$ to the Hamiltonian
side, i.e., expressing it in the variables $(\mathbf{m}, A)$, it
follows that $\{h, j_m\} = 0$. It is important to note that this
$U(1)$-action with momentum map $j_m$ is expressed in the variables
$(\mathbf{m}, A)$ as: $(\mathbf{m}, A) \mapsto \left(\mathbf{m}, 
e^{ {\rm i}\varphi}A \right)$.

Now, consider a second circle action on 
$\mathfrak{so}(3) ^\ast\times \operatorname{Orb}(A_0)$ given by 
\[
(\mathbf{m}, A) \mapsto (\rho(\varphi) \mathbf{m}, \rho(\varphi)A\rho(\varphi)^{-1}),\quad \text{where} \quad  
\rho(\varphi) : = \exp (\varphi \widehat{\mathbf{e}}_3).
\]
 This is
the coadjoint action of a subgroup of $SO(3) \,\circledS\, 
\mathfrak{gl}(3, \mathbb{C})$ and hence it is Poisson. It admits
a momentum map which is $\mathbf{J}_3$. The Hamiltonian $h$ given 
by \eqref{Hamiltonian_A_second} is invariant under this action and
so we conclude that $\{h, \mathbf{J}_3\}=0$. The action 
$(\mathbf{m}, A) \mapsto (\rho(\varphi) \mathbf{m}, \rho(\varphi)A\rho(\varphi)^{-1})$ is induced via the cotangent bundle
version of Theorem \ref{EP_reduction_corollary}
by the cotangent lift of the
action $A \mapsto \rho(\varphi)A\rho(\varphi)^{-1}$. This
action on configuration space commutes 
with the previously considered circle action $A \mapsto e^{{\rm i} \varphi}A$. Therefore, the associate momentum maps commute, i.e., 
$\{j_m, \mathbf{J}_3\}=0$. Concluding, we have
$\{h, j_m\} = 0$, $\{h, \mathbf{J}_3\}=0$, $\{j_m, \mathbf{J}_3\}=0$.

Finally, we prove the functional independence of the three integrals
$h, j_m, \mathbf{J}_3$. Instead of showing that their differentials
are linearly independent away from a subset of measure zero in
$\mathfrak{so}(3)^* \times \operatorname{Orb}(A_0)$, we will show that
the Hamiltonian vector fields generated by these integrals are
independent on such a set. Since $j_m$  and $\mathbf{J}_3$ are  
momentum maps, their Hamiltonian vector fields relative to the 
Lie-Poisson bracket \eqref{Lie_Poisson_bracket_A_second} coincide with 
the infinitesimal generator vector fields of the corresponding 
$U(1)$-actions. These vector fields are hence $(\mathbf{m}, A) \mapsto
(\mathbf{m}, A; {\bf 0}, {\rm i}A)$ and $(\mathbf{m}, A) \mapsto
\left(\mathbf{m}, A; \mathbf{e}_3\times \mathbf{m}, 
\left[\widehat{\mathbf{e}}_3, A\right]\right)$.

Now we compute the Hamiltonian vector field for $h$ given by
\eqref{Hamiltonian_A_second}. We have $\delta h/ \delta \mathbf{m} = 
\frac{1}{2} \mathbf{J}(A)^{-1} \mathbf{m}$. A direct computation shows
that 
\[
\frac{\delta h}{\delta A} = 2\left[\widehat{\mathbf{w}}, \left[
\widehat{\mathbf{w}}, A \right]\Gamma \right], \quad \text{where}
\quad \mathbf{w}: = \frac{1}{2}\mathbf{J}(A) ^{-1} \mathbf{m}.
\]
Therefore, from \eqref{LP_A_second_Ham_equ}, we obtain the
expression of the Hamiltonian vector field defined by $h$, namely,
\[
(\mathbf{m}, A) \longmapsto X_h(\mathbf{m}, A) = 
\left(\mathbf{m}, A; \mathbf{w} \times \mathbf{m} - 
4 \operatorname{Re}\overrightarrow{\big[\left[\Gamma \left[
A, \widehat{\mathbf{w}}\right],\widehat{\mathbf{w}} \right],  A \big]
\phantom{\Big]}\!\!},
\left[\widehat{\mathbf{w}}, A \right]\right) .
\]
We need to show that
\[
\left\{
\begin{aligned}
&\alpha_1 \left( \mathbf{w}\times \mathbf{m} - 4\operatorname{Re}
\overrightarrow{X^{\phantom{4}}\!\!}\right) + \alpha_3 \mathbf{e}_3 \times\mathbf{m} = 0,\\
&\alpha_1\left[\widehat{\mathbf{w}}, A \right] + \alpha_2{\rm i}A
+ \alpha_3\left[\widehat{\mathbf{e}}_3, A \right] = 0,
\end{aligned}
\right.
\]
where $X: = \big[\left[\Gamma \left[
A, \widehat{\mathbf{w}}\right],\widehat{\mathbf{w}} \right],  A \big]$,
implies that $\alpha_1 = \alpha _2 = \alpha_3 = 0$ for generic 
$(\mathbf{m}, A)$. Taking the dot product of the first equation 
with $\mathbf{m}$ yields $\alpha_1\operatorname{Re}
\overrightarrow{X^{\phantom{4}}\!\!} \cdot \mathbf{m} = 0$. It is
easy to find points $(\mathbf{m}, A)$ for which 
$\operatorname{Re}\overrightarrow{X^{\phantom{4}}\!\!} \cdot 
\mathbf{m} \neq 0$. Since this expression is polynomial in
$\mathbf{w}$ and $A$ and since it does not vanish identically,
its set of zeros is of measure zero in $\mathfrak{so}(3)^*
\times \operatorname{Orb}(A_0)$. This shows that for a set of
measure zero on this phase space, $\alpha_1 =0$. Choosing
$\mathbf{m}$ not collinear with $\mathbf{e}_3$, implies that
$\alpha_3 = 0$. Now, for any $A \in \operatorname{Orb}(A_0) 
\neq 0$, we get $\alpha_2 =0$. We have proved the following 
result.

\begin{theorem}
The three functions $h, j_m ,  \mathbf{J}_3$ form a completely 
integrable system on the six dimensional coadjoint orbit
$\mathfrak{so}(3)^* \times \operatorname{Orb}(A_0)$.
\end{theorem}

\subsubsection{The B-phase}

We consider the element
\[
I _3  = \left(
\begin{array}{ccc}
1&0&0\\
0&1&0\\
0&0&1
\end{array} 
\right)\in \mathfrak{gl}(3, \mathbb{C}).
\]

\begin{proposition} {\rm (i)} The isotropy subgroup of $I _3 $ is 
\[
H=\{(1,R , R)\mid R \in SO(3)\} \subset G=U(1) \times SO(3)_L \times SO(3)_R.
\]
{\rm (ii)} We have the diffeomorphism
\begin{equation}\label{diffeo_B1} 
G/H \ni [e^{{\rm i}\varphi }, R _1, R _2 ]_H\longmapsto 
(e^{{\rm i}\varphi }, R _1 R _2^{-1}) \in U(1) \times SO(3).
\end{equation} 
{\rm (iii)} We have the diffeomorphism
\begin{equation}\label{diffeo_B2} 
U(1) \times SO(3)\ni (e^{{\rm i}\varphi },R) \longmapsto 
e^{{\rm i}\varphi }R\in \operatorname{Orb}(I_3).
\end{equation} 
\end{proposition} 
\textbf{Proof.} {\rm (i)} By taking the imaginary part of the equality $e^{ {\rm i} \varphi }R _1 R _2 ^{-1} =I_3$, we obtain that    $ \varphi\in \{0, \pi\}$. So we get $\pm R _1 R _2 =I_3$. Taking the determinant shows that the minus sign is impossible, therefore the result follows.\\
{\rm (ii)} The result follows from a direct verification.\\
{\rm (iii)} From (i), it follows that we have the diffeomorphism $[e^{{\rm i}\varphi }, R _1, R _2 ]_H \in G/H \mapsto e^{ {\rm i} \varphi }R _1 R _2 ^{-1} \in  \operatorname{Orb}(I_3)$. The result then follows by composing with the diffeomorphism obtained in (ii). $ \qquad\blacksquare$  

\medskip

We observe that the subgroup $\widetilde{G}: = U(1) \times SO(3)_L \subset G$ acts transitively on the orbit $\operatorname{Orb}(I_3)$, 
see \eqref{diffeo_B1}, \eqref{diffeo_B2}. The isotropy subgroup 
of $I_3 $ is $\widetilde{G}_{I _3 }=H \cap \widetilde{G}=
\{1, I _3 , I _3 \}$. We recover the fact that the orbit 
$\operatorname{Orb}(I_3) \subset \mathfrak{gl}(3, \mathbb{C})$ 
is diffeomorphic to $U(1) \times SO(3)$.
\medskip

 As a consequence, we apply Theorem \ref{EP_reduction_corollary} with this description of the orbit and hence the Lie algebra one has to consider is $\mathbb{R}   \times \mathfrak{so}(3)_L$.
On this orbit $M$ we consider the Lagrangian density given by the gradient part only, i.e. $\mathcal{L}(A, \nabla A) = F_{grad}(A, \nabla A)$. 

The reduced velocity $\xi(z)  = \partial _z  g g^{-1} $ of the general theory (see 
\eqref{reduced_section_formula_for_g}) is given 
here by $ \xi = (v,\mathbf{w}) : X  \rightarrow  \mathbb{R}  \times \mathfrak{so}(3)_L $, where $v= \partial _z \varphi$, and $\mathbf{w} =  R^{-1}(\partial_z R )$. The Euler-Poincar\'e Lagrangian
\[
l=l(\xi, m): \mathbb{R}  \times \mathfrak{so}(3)_L  \times M \rightarrow 
\mathbb{R}
\]
of the general theory given in  \eqref{EP_Lagrangian}, is computed in this case to be
\begin{align*} 
l(v,\mathbf{w},A) &=F_{grad}(A, \partial_z A)=
\left\langle \! \left\langle {\rm i}v A+A \widehat{\mathbf{w}},
{\rm i}v A+A \widehat{\mathbf{w}} \right\rangle \! \right\rangle  \\
&=v ^2 \left\langle \! \left\langle A , A \right\rangle \! \right\rangle + \mathbf{w}^\mathsf{T} I(A) \mathbf{w} +2 v \left\langle \! \left\langle A {\rm i}, A \widehat{ \mathbf{w} } \right\rangle \! \right\rangle \\
&= 2 \gamma _1 \| \mathbf{w}\| ^2 +( \gamma _2 + \gamma _3 ) (w _1 ^2 + w _2 ^2 )+( 3 \gamma _1 + \gamma _2 + \gamma _3 ) v ^2 ,
\end{align*} 
where $I_{ab}(A)= \left\langle\!\left\langle A \widehat{\mathbf{e} }_a, A \widehat{\mathbf{e} } _b \right\rangle \! \right\rangle $ and we used $A= e^{{\rm i}\varphi }R$ so that $A ^\ast A=I_3$.

The Euler-Poincar\'e equations \eqref{EP_equations} are
\begin{equation}\label{EP_B} 
\partial _z \frac{\delta l}{\delta v } = \mathbf{J}_1   \left( \frac{\delta l}{\delta A} \right), \qquad \partial _z \frac{\delta l}{\delta \mathbf{w} } - 
\operatorname{ad}^*_ \mathbf{w}  \frac{\delta l}{\delta \mathbf{w} }= 
\mathbf{J} _2  \left( \frac{\delta l}{\delta A} \right),
\end{equation} 
where $ \mathbf{J} _1:T^*M \rightarrow \mathbb{R} $ is the momentum 
map of the $U(1)$ action $A \mapsto e^{{\rm i}\psi  }A$, and 
$ \mathbf{J} _2 : T^*M \rightarrow \mathfrak{so}(3)_L^\ast $ is 
the momentum map of the $SO(3)$ action $A \rightarrow QA$. They 
have the expressions
\begin{equation}\label{momentum_map_B}
\mathbf{J}_1 ( \alpha _A)= \operatorname{Re} 
\operatorname{Tr}(\alpha_A^\ast{\rm i}A) \quad\text{and}\quad 
\mathbf{J}_2 (\alpha_A )=2\overrightarrow{\operatorname{Re}
\left(\alpha_A^\ast A\right)^{\phantom{4}}\!\!}.
\end{equation}
The second equation in \eqref{EP_equations} is given, in this case, by
\[
\partial _z A= {\rm i}vA+ A\widehat{\mathbf{w}}.
\]

Since
\[
\frac{\delta l}{\delta v}=2( 3 \gamma _1 + \gamma _2 + \gamma _3 ) v, \quad \frac{\delta l}{\delta \mathbf{w} }=4 \gamma _1 \mathbf{w} +2( \gamma _2+ \gamma _3 )(w _1 , w _2 , 0 )^\mathsf{T}= \mathbf{J} \mathbf{w} , \quad \frac{\delta l}{\delta A}=0,
\]
where $ \mathbf{J} = \operatorname{diag}(4 \gamma _1 + 2 \gamma _2 + 2 \gamma _3 ,4 \gamma _1 + 2 \gamma _2 + 2 \gamma _3,4 \gamma _1 )$, we get
\begin{equation}
\label{EP_B_phase}
\partial _z v= 0, \qquad \partial _z\mathbf{J} \mathbf{w}  + \mathbf{w} \times  \mathbf{J} \mathbf{w} =0, \qquad \partial _z A= ivA+ A\widehat{ \mathbf{w} }.
\end{equation} 
The equation for $ \mathbf{w} $ is the Euler equation for a free symmetric rigid body with moment of inertia $ \mathbf{J} $ whose solutions are well-known (see e.g., \cite[\S5.7]{Lawden1989}). In this special case, all equations decouple and \eqref{EP_B_phase} can 
be solved explicitly.

\begin{remark}{\rm 
An interesting question is the transition from one phase to
another. In the case of liquid crystals, this means the transition from
biaxial nematics to uniaxial nematics which leads to the deformation
of the orbit $SO(3)/(\mathbb{Z}_2 \times \mathbb{Z}_2)$ to 
$SO(3)/\mathbb{Z}_2$. We do not consider this problem in the present paper. \quad $\blacklozenge$ }
\end{remark}

\subsection{Neutron stars}

The classification of the phases in neutron stars parallels the procedure
used in the classification of the phases for superfluid liquid Helium
${}^3$He, taking into account the differences associated to the order 
parameter of these two problems. For details, see \cite{MoSa2011}. We
shall use the notations of this paper below and consider several examples
of orbits.
\medskip

Define $S^2(\mathbb{C}^3)_0: = \left\{A \in \mathfrak{gl}(3,\mathbb{C}) 
\mid A^\mathsf{T} = A, \; \operatorname{Tr}(A)=0 \right\}$. The direct
product Lie group $U(1) \times SO(3)$ acts on $S^2(\mathbb{C}^3)_0$ by
\begin{equation}
\label{neutron_star_action}
A \longmapsto e^{{\rm i}\varphi} RAR ^{-1},
\end{equation}
where $e^{{\rm i}\varphi} \in U(1)$ and $R\in SO(3)$. The possible 
non-trivial orbits of this group action are given in \cite{MoSa2011}
and denoted by $\Omega_i$, where $i=1, \ldots, 10$.
We shall consider below only the most interesting examples $i=1,4, 6,8$.

\subsubsection{$\Omega_1$-phase}

Consider the orbit of the group $U(1) \times SO(3)$ through 
\[
A_0 = 
\begin{pmatrix}
1& 0&\;\;0\\
0& 1& \;\;0\\
0&0&-2
\end{pmatrix} \in S^2(\mathbb{C}^3)_0.
\]
A direct verification shows that the isotropy subgroup of $A_0$
for the action \eqref{neutron_star_action} is the one-dimensional
subgroup $\left\{(1, \rho(\varphi)J_\pm) \mid \varphi \in  \mathbb{R}
\right\} \subset U(1) \times SO(3)$, where 
$\rho(\varphi) = \exp(\varphi\widehat{\mathbf{e}}_3)$, $J_\pm$ 
is defined in \eqref{J_tilde_J}, and $\rho(\varphi)J_\pm = J_\pm\rho(\varphi)$. 
So, the $(U(1) \times SO(3))$-orbit 
though $A_0$ is diffeomorphic to $U(1) \times \left[SO(3)/(\mathbb{Z}_2
\times SO(2)) \right]$ which is the $\Omega_1$-orbit in
the classification \cite{MoSa2011}. Explicitly, this 
$U(1) \times SO(3)$-orbit is
\[
M=\operatorname{Orb}(A_0) = \left\{e^{{\rm i}\varphi}\left(
\mathbf{u}\mathbf{u}^\mathsf{T} + \mathbf{v}\mathbf{v}^\mathsf{T} -
2\mathbf{w}\mathbf{w}^\mathsf{T}\right) \mid \mathbf{u}, \mathbf{v} \in 
S^2, \mathbf{u}\cdot \mathbf{v}= 0, \mathbf{w} = \mathbf{u}\times \mathbf{v} \right\}.
\]
Note that, in accordance with the isotropy subgroup described above, 
the identities characterizing $M=\operatorname{Orb}(A_0)$ are unchanged
if we replace \\
$\bullet$ $\mathbf{u}$ by $\mathbf{u} \cos \theta + \mathbf{v}\sin\theta$ 
and $\mathbf{v}$ by $- \mathbf{u} \sin\theta + 
\mathbf{v} \cos\theta$ (this is the $SO(2)$-action)\\
$\bullet$ $\mathbf{v}$ by $- \mathbf{v}$ and hence $\mathbf{w}$ by
$-\mathbf{w}$ (this is the $J_-$-action).

Note that since the isotropy subgroup is Abelian, its coadjoint 
orbits are points. Therefore, by Theorem \ref{symplectic_leaves}, 
the canonical cotangent bundle $T^*M$ is symplectically diffeomorphic
to the coadjoint orbit $\mathcal{O}_{(p, \mathbf{m}, A_0)}$ endowed,
as usual, with the minus orbit symplectic structure, $(p, \mathbf{m})
\in \left(\mathbb{R} \times \mathfrak{so}(3) \right)^*$.
\medskip

To find the Lagrangian associated to this phase, we note that the
reduced velocity $\xi= \partial _z g g^{-1}$ of the general theory (see 
\eqref{reduced_section_formula_for_g}) is given 
here by $ \xi = (v,\mathbf{w} ) : \mathbb{R}  \rightarrow  \mathbb{R}
\times \mathfrak{so}(3)$, where $v = \partial_z \varphi$ and 
$\mathbf{w}=(\partial_z R) R^{-1}$. The infinitesimal generator of
the group action, i.e., the second formula in \eqref{EP_equations}, 
has in this case the expression
\[
\partial _z A= {\rm i}v A+ \left[\widehat{\mathbf{w}}, A \right].
\] 
This formula and \eqref{Fgrad_1D} show that the Euler-Poincar\'e Lagrangian \eqref{EP_Lagrangian} of the general theory becomes in
this case the function $l=l(\xi, m):  \mathbb{R}\times \mathfrak{so}(3)  
\times M \rightarrow \mathbb{R}$ given by
\begin{align}
\label{Lagrangian_omega_1} 
l(v ,\mathbf{w} ,A) &= \left\langle \! \left\langle
{\rm i}v A+ \left[\widehat{\mathbf{w}}, A \right], 
{\rm i}v A+ \left[\widehat{\mathbf{w}}, A \right]
\right\rangle \! \right\rangle \nonumber \\
&=v^2 \left\langle \! \left\langle A, A \right\rangle \! \right\rangle
+ 2v \left\langle \! \left\langle {\rm i} A, 
\left[\widehat{\mathbf{w}}, A \right]
\right\rangle \! \right\rangle + 
\mathbf{w}^\mathsf{T} \mathbf{J}(A) \mathbf{w}
\end{align} 
(see \eqref{inner_product} and \eqref{J_A}). Thus,
the Euler-Poincar\'e equations \eqref{EP_equations} read
\begin{equation}
\label{EP_Omega_1}
\begin{aligned}
&\partial _z \frac{\delta l}{\delta v} = 
\operatorname{Re}\operatorname{Tr}\left(\frac{\delta l}{\delta A}^* A
{\rm i}\right), \qquad 
\partial_z \frac{\delta l}{\delta \mathbf{w}} +
\frac{\delta l}{\delta \mathbf{w}} \times \mathbf{w} = 
2\overrightarrow{\operatorname{Re}\left[\frac{\delta l}{\delta A}^*,  A\right]^{\phantom{4}}\!\! } ,\\ 
&\partial _z A= {\rm i}v A+ \left[\widehat{\mathbf{w}}, A \right],
\end{aligned}
\end{equation}
and where we have
\begin{align*}
\frac{\delta l}{\delta v}&=2v \left\langle \! \left\langle A, A
\right\rangle \! \right\rangle + 2\left\langle \! \left\langle 
{\rm i} A, \left[\widehat{\mathbf{w}}, A \right]
\right\rangle \! \right\rangle ,  \qquad 
\frac{\delta l}{\delta \mathbf{w}}=2 \mathbf{J}(A) \mathbf{w} -
4v \overrightarrow{\operatorname{Re}\left[{\rm i} \Gamma A^*,  A\right]^{\phantom{4}}\!\! }, \\
\frac{\delta l}{\delta A}&= \left[\left[\widehat{\mathbf{w} },A\right] 
\Gamma ,\widehat{\mathbf{w} }\right] + 
\left[\left[\widehat{\mathbf{w} },A\right] 
\Gamma ,\widehat{\mathbf{w} }\right]^\mathsf{T} + X,
\end{align*}
where $X$ is the traceless symmetric part of $v^2 A\Gamma - 2v {\rm i}
\left[\widehat{\mathbf{w}}, A \right] \Gamma - 2v {\rm i}
\left[\widehat{\mathbf{w}}, A\Gamma \right]$.

The Hamiltonian associated to $l$ on the six dimensional coadjoint orbit 
$\mathcal{O}_{(p, \mathbf{m}, A_0)}$ generates an integrable system.
The integrals in involution are $h, j_m, \mathbf{J}_3$ (see the text 
after \eqref{LP_A_second_Ham_equ} for the expressions of $j_m$ and
$\mathbf{J}_3$). To show independence on a dense open subset of phase
space, we shall work on the Hamiltonian side. As in Subsection 
\ref{sec:A_phase_second}, the two circle actions
\[
(p, \mathbf{m}, A) \mapsto \left(p, \mathbf{m}, e^{{\rm i} \varphi}A\right), \qquad  (p, \mathbf{m}, A) \mapsto (p, \rho(\varphi) \mathbf{m}, 
\rho(\varphi)A\rho(\varphi)^{-1}),
\]
where $\rho(\varphi) : = \exp (\varphi \widehat{\mathbf{e}}_3)$, 
generate the momentum maps $j_m$ and $\mathbf{J}_3$, respectively.
The infinitesimal generator vector fields of these actions coincide
with the Hamiltonian vector fields of $j_m$ and $\mathbf{J}_3$,
respectively, i.e., they are
\[
X_{j_m}(p, \mathbf{m}, A) =
(0,{\bf 0}, {\rm i}A), \qquad X_{\mathbf{J}_3}(p, \mathbf{m}, A) 
= \left(0, \mathbf{e}_3\times \mathbf{m}, 
\left[\widehat{\mathbf{e}}_3, A\right]\right).
\] 
Finally the Hamiltonian vector field is obtained
by using \eqref{EP_Omega_1} and the relations $\frac{\delta h}{\delta A}
= - \frac{\delta l}{\delta A}$, $\mathbf{w} =\frac{\delta h}{\delta \mathbf{m}}$, $v = \frac{ \delta h}{ \delta v}$. A direct computation
yields
\begin{align*}
X_h(p, \mathbf{m}, A) = \left(
-\operatorname{ReTr} \left(\frac{\delta h}{\delta A}^* {\rm i}A\right),\, \frac{\delta h}{\delta \mathbf{m}} \times \mathbf{m} - 
2\overrightarrow{\operatorname{Re}\left[\frac{\delta h}{\delta A}^*,  A\right]^{\phantom{4}}\!\! },\, {\rm i}\frac{\delta h}{\delta p} A + 
\left[\widehat{\frac{\delta h}{\delta \mathbf{m}}}, A \right] 
\right).
\end{align*}
Writing $\alpha_1 X_h(p, \mathbf{m}, A) + \alpha_2 X_{\mathbf{J}_3}(p, \mathbf{m}, A) + \alpha_3 X_{j_m}(p, \mathbf{m}, A) = 0$, implies that
$\alpha_1 = \alpha_2 = \alpha_3 = 0$ on an open dense set of
points $(p, \mathbf{m}, A)$ in the coadjoint orbit.

\begin{theorem}
The three functions $h, j_m ,  \mathbf{J}_3$ form a completely 
integrable system on all six dimensional coadjoint orbits
$\mathcal{O}_{(p, \mathbf{m}, A_0)}$.
\end{theorem}

\subsubsection{$\Omega_4$-phase}
The $\Omega_4$-phase corresponds to 
\[
A_0 = 
\begin{pmatrix}
1&0&0\\
1&\omega&0\\
0&0&\omega^2
\end{pmatrix} \in S^2(\mathbb{C}^3)_0,
\]
where $\omega^3 = 1$. The isotropy subgroup is
\begin{align*}
\big(U(1) \times SO(3)\big)_{A_0} &= 
\left.\left\{\left(1, 
\begin{pmatrix}
\epsilon_1&0&0\\
0&\epsilon_2&0\\
0&0&\epsilon_3
\end{pmatrix} \right), \left(\omega, 
\begin{pmatrix}
0&0&\epsilon_1\\
\epsilon_2&0&0\\
0&\epsilon_3&0
\end{pmatrix} \right), \left( \omega^2, 
\begin{pmatrix}
0&\epsilon_1&0\\
0&0&\epsilon_2\\
\epsilon_3&0&0
\end{pmatrix} \right) \, \right| \right.\\
& \left.\phantom{
\begin{pmatrix}
0&0&\epsilon_1\\
\epsilon_2&0&0\\
0&\epsilon_3&0
\end{pmatrix}  }
\epsilon_i = \pm 1,\; i=1,2,3,\;
\epsilon_1 \epsilon_2 \epsilon_3 = 1 \right\}
\end{align*}
and it is isomorphic to the tetrahedral group (the 12 elements alternating 
group $\mathfrak{A}_4$ on 4 letters).
The $(U(1) \times SO(3))$-orbit $\operatorname{Orb}(A_0) = M$ is 
hence diffeomorphic to $\left[U(1) \times SO(3)\right]/\mathfrak{A}_4$,
which is the $\Omega_4$-coadjoint orbit in
the classification \cite{MoSa2011}.
\medskip

To compute the texture equations associated to this phase, we note that
the reduced velocity $\xi = \partial _z  g g^{-1}$ of the general theory (see 
\eqref{reduced_section_formula_for_g}) is given 
here by $ \xi = (v,\mathbf{w} ) : \mathbb{R}  \rightarrow  \mathbb{R}
\times \mathfrak{so}(3)$, where $v = \partial_z \varphi$ and 
$\mathbf{w}=(\partial_z R) R^{-1}$. The second formula in \eqref{EP_equations} 
(the infinitesimal generator of the action) is given here by
\[
\partial _z A= {\rm i}v A+ \left[\widehat{\mathbf{w}}, A \right].
\] 
Using this expression and formula \eqref{Fgrad_1D}, the Euler-Poincar\'e Lagrangian
\[
l=l(\xi, m):  \mathbb{R}\times \mathfrak{so}(3)  \times M \rightarrow 
\mathbb{R}
\]
of the general theory given in  \eqref{EP_Lagrangian}, is computed in this case to be
\begin{align}
\label{Lagrangian_omega_4} 
l(v ,\mathbf{w} ,A) &= \operatorname{Re} \operatorname{Tr}( \Gamma \partial _z A ^\ast \partial _z A) \nonumber \\
&=  \mathbf{w}^\mathsf{T} \mathbf{J}(A) \mathbf{w} + 
2v \left\langle \! \left\langle {\rm i} A, 
\left[\widehat{\mathbf{w}}, A \right]
\right\rangle \! \right\rangle +(3\gamma_1 + \gamma_2 + \gamma_3)v^2,
\end{align} 
(see \eqref{inner_product} and \eqref{J_A}). Thus,
the Euler-Poincar\'e equations \eqref{EP_equations} read
\begin{equation*}
\partial _z \frac{\delta l}{\delta v} = 
\operatorname{Re}\operatorname{Tr}\left(\frac{\delta l}{\delta A}^* A
{\rm i}\right), \quad 
\partial_z \frac{\delta l}{\delta \mathbf{w}} +
\frac{\delta l}{\delta \mathbf{w}} \times \mathbf{w} = 
2\overrightarrow{\operatorname{Re}\left[\frac{\delta l}{\delta A}^*,  A\right]^{\phantom{4}}\!\! } ,\quad 
\partial _z A= {\rm i}v A+ \left[\widehat{\mathbf{w}}, A \right],
\end{equation*}
 where we have
\begin{align*}
&\frac{\delta l}{\delta v}=2(3\gamma_1+ \gamma_2 + \gamma_3)v
+ \left\langle \! \left\langle {\rm i} A, 
\left[\widehat{\mathbf{w}}, A \right]
\right\rangle \! \right\rangle,  \qquad 
\frac{\delta l}{\delta \mathbf{w}}=2 \mathbf{J}(A) \mathbf{w}
- 4v \overrightarrow{\operatorname{Re}\left[{\rm i} \Gamma A^*,  A\right]^{\phantom{4}}\!\! }, \\
& \frac{\delta l}{\delta A}= \left[\left[\widehat{\mathbf{w} },A\right] 
\Gamma ,\widehat{\mathbf{w} }\right] + 
\left[\left[\widehat{\mathbf{w} },A\right] 
\Gamma ,\widehat{\mathbf{w} }\right]^\mathsf{T} + X,
\end{align*}
where $X$ is the traceless symmetric part of $v^2 A\Gamma - 2v {\rm i}
\left[\widehat{\mathbf{w}}, A \right] \Gamma - 2v {\rm i}
\left[\widehat{\mathbf{w}}, A\Gamma \right]$.

Passing to the Hamiltonian formulation, a direct computation using
\eqref{red_Poisson}  gives the Lie-Poisson bracket
\begin{align}
\label{Lie_Poisson_bracket_omega_4}
\{f,h\}(p,\mathbf{m} , A)= &\mathbf{m} \cdot \left(
\frac{\delta f}{\delta \mathbf{m}}\times 
\frac{\delta h}{\delta \mathbf{m}} \right)+ 
\left\langle \frac{\delta h}{\delta A}, 
\left[ A,\widehat{\frac{\delta f}{\delta \mathbf{m}}}\right] 
-{\rm i}\frac{\delta f}{\delta p} A \right\rangle \nonumber \\
& \qquad - \left\langle \frac{\delta f}{\delta A}, 
\left[A,\widehat{\frac{\delta h}{\delta \mathbf{m}}}\right] 
- {\rm i}\frac{\delta h}{\delta p} A \right\rangle\,,
\end{align}
where $p :=\frac{\delta l}{\delta v}=2(3\gamma_1+ \gamma_2 + \gamma_3)v$
and $\mathbf{m}:= \frac{\delta l}{\delta \mathbf{w}}=
2 \mathbf{J}(A) \mathbf{w}$. Thus, the equations 
$\partial_z f = \{f, h\}$ for any $f$ are 
\begin{equation}
\label{LP_omega_4_Ham_equ}
\begin{aligned}
\partial_z p &= - \operatorname{ReTr} \left(\frac{\delta h}{\delta A}^*
A {\rm i}\right), \qquad
\partial_z\mathbf{m}+ \mathbf{m} \times 
\frac{\delta h}{\delta \mathbf{m}} = 
-2\operatorname{Re}\overrightarrow{\left[\frac{\delta h}{\delta A} ^\ast, A \right]^{\phantom{4}}\!\!} \\
\partial_z A &= {\rm i} \frac{\delta h}{\delta p} A +\left[\widehat{\frac{\delta h}{\delta \mathbf{m}}},A\right].
\end{aligned}
\end{equation}

The Hamiltonian system generated by $h$ has three integrals of motion
in involution: $h, \mathbf{J}_3, j_m = p$. The integrals $\mathbf{J}_3$
and $j_m$ commute with $h$ because on the Lagrangian side they are
momentum maps and hence constant on the solutions of the Euler-Lagrange
equations. The maps $\mathbf{J}_3$ and $j_m$ commute because they
are momentum maps of two commuting circle actions (this is the
same argument as given in \S\ref{sec:A_phase_second}). The fact that
$p = j_m$ is a direct computation replacing $\partial_z A = {\rm i}vA + 
\left[\widehat{\mathbf{w}}, A \right]$ in the defining formula
\eqref{j_m_A_first} for $j_m$ and using $A = 
e^{{\rm i}\varphi}RA_0R^{-1}$, $A_0A_0^* = I_3$. For the complete
integrability of this system, one more integral is needed.

\subsubsection{$\Omega_6$-phase}
 Consider the
orbit of this group through 
\[
A_0 = 
\begin{pmatrix}
1& \;\;{\rm i}&0\\
{\rm i}& -1& 0\\
0&\;\;0&0
\end{pmatrix} \in S^2(\mathbb{C}^3)_0.
\]
We note that 
\begin{equation}
\label{rotation_identity}
e^{{\rm i}\varphi} A_0 = \rho\left(-\frac{\varphi}{2}\right)A_0 
\rho\left(\frac{\varphi}{2}\right), 
\end{equation}
where $\rho(\varphi)= \operatorname{exp}(\varphi\widehat{\mathbf{e}}_3)$.
A direct verification shows that the isotropy subgroup of $SO(3)$
for the action \eqref{neutron_star_action} is $\{\tilde{J}_\pm\} 
\cong \mathbb{Z}_2$ (see \eqref{J_tilde_J}). As a consequence, using
\eqref{rotation_identity}, it follows that the isotropy subgroup
of the action \eqref{neutron_star_action} equals  
\[
\left(U(1) \times SO(3)\right)_{A_0} = 
\left.\left\{ \left( e^{{\rm i}\varphi}, \tilde{J}_\pm \rho\left(
\frac{ \varphi}{2}\right) \right) \,\right| \, 
e^{{\rm i}\varphi} \in U(1)
\right\} \cong U(1) \times \mathbb{Z}_2.
\]

The group $\{\tilde{J}_\pm\} \cong \mathbb{Z}_2=\{\pm 1\}$ acts on 
$SO(3)$ by $R \mapsto R\tilde{J}_\pm$ and $\left(U(1) \times SO(3)\right)_{A_0}$ on $U(1) \times SO(3)$ by 
$\left(e^{{\rm i} \psi}, R \right) \mapsto \left(e^{{\rm i} \psi}, 
R \right)\left(e^{{\rm i}\varphi}, \tilde{J}_\pm \rho\left(
\frac{\varphi}{2}\right) \right)$.
Thus 
\[
\left(U(1) \times SO(3)\right)/\left(U(1) \times SO(3)\right)_{A_0}
\ni \left[e^{{\rm i} \psi}, R \right] \longmapsto \left[R \rho\left(-
\frac{\varphi}{2}\right)\right] \in SO(3)/\mathbb{Z}_2
\]
is a diffeomorphism which shows that the orbit 
$\left(U(1) \times SO(3)\right)\cdot A_0$ is 
diffeomorphic to  $SO(3)/\mathbb{Z}_2$. We have found the orbit
of type $\Omega_6$ in the classification given in \cite{MoSa2011}. 

It is easy to see that 
\[
R^\mathsf{T} A_0 R = (R_1 + {\rm i}R_2)(R_1 + {\rm i}R_2)^\mathsf{T},
\]
where $R_j$ is the $j^{\rm th}$ column of $R \in SO(3)$, $j=1,2$. In view
of the previous considerations this shows that 
\begin{equation}
\label{orb_a_0_neutron}
\operatorname{Orb}(A_0) = \left\{(\mathbf{x} + {\rm i}\mathbf{y})
(\mathbf{x} + {\rm i}\mathbf{y})^\mathsf{T}\mid \mathbf{x}, \mathbf{y}\in 
S^2,\, \mathbf{x}\cdot \mathbf{y}= 0\right\} \cong SO(3)/\mathbb{Z}_2.
\end{equation}

\paragraph{Lagrangian formulation.} Note that the setting is very similar
to the setup for superfluid liquid Helium ${}^3$He in the second regime
of phase A, since the action turns out the be the same, namely,  $A
\mapsto RAR^{-1}$. The orbit is, however, different because the matrices
$A_0$ in these two cases do not lie on the same $SO(3)$-orbit.

The Lagrangian used in this case (see \cite{MoSa2011}) is given by
\eqref{Fgrad_1D}. Therefore, the Euler-Poincar\'e Lagrangian takes
the same form as in the second regime of the A-phase for superfluid liquid Helium ${}^3$He, namely,  
\[
l(\mathbf{w} ,A) = \sum_{a,b=1}^3 J_{ab}(A) w _a w _b = \mathbf{w} ^\mathsf{T}\mathbf{J} (A)\mathbf{w},
\]
where $\mathbf{J}(A)$ is given in \eqref{J_A}.
The Euler-Poincar\'e equations are thus \eqref{equ_A_phase_second}. The
functional derivatives are in this case
\[
\frac{\delta l}{\delta \mathbf{w} }=2 \mathbf{J}(A) \mathbf{\mathbf{w} }, \quad \frac{\delta l}{\delta A}= [[\widehat{\mathbf{w} },A] \Gamma ,\widehat{\mathbf{w}}] + [[\widehat{\mathbf{w} },A] \Gamma ,\widehat{\mathbf{w}}]^\mathsf{T} \in S^2(\mathbb{C}^3)_0.
\]
\paragraph{Hamiltonian formulation.} All formulas, with the changes
noted above, are identical to the ones in the second regime
of phase A  for superfluid liquid Helium ${}^3$He, i.e., \eqref{Hamiltonian_A_second}
and \eqref{Lie_Poisson_bracket_A_second} hold. The same considerations
about the complete integrability of the equations given at the
end of \S\ref{sec:A_phase_second} hold because the action given by 
multiplication with $e^{{\rm i}\varphi}$ preserves the orbit $\operatorname{Orb}(A_0)$ given by \eqref{orb_a_0_neutron} in view of
\eqref{rotation_identity}. We get the following result.

\begin{theorem}
The three functions $h, j_m ,  \mathbf{J}_3$ form a completely 
integrable system on the six dimensional coadjoint orbit
$\Omega_6$.
\end{theorem}

\subsubsection{$\Omega_8$-phase}
The $\Omega_8$-phase corresponds to 
\[
A_0^\pm = 
\begin{pmatrix}
\;\;0&1&\pm {\rm i}\\
\;\;1&0&\;\;0\\
\pm {\rm i}&0&\;\;0
\end{pmatrix} \in S^2(\mathbb{C}^3)_0.
\]
Note that $e^{\mp{\rm i} \varphi} A_0^\pm = \rho(\varphi)A_0^\pm
\rho(-\varphi)$, where $\rho(\varphi) : = 
\exp(\varphi\widehat{\mathbf{e}}_1)$. Hence $SO(3)$ acts transitively
on the orbit through $A_0$, i.e., $\operatorname{Orb}(A_0) = 
\{RA_0 R ^{-1} \mid R \in SO(3)\}$. A direct computation shows that
the isotropy group $SO(3)_{A_0} = \{I_3\}$ and hence 
$\operatorname{Orb}(A_0)$ is diffeomorphic to $SO(3)$.

It is easy to verify that the isotropy subgroup of the original action 
equals
\[
(U(1) \times SO(3))_{A_0^\pm} = \left\{\left(e^{{\rm i} \varphi}, 
\rho(\pm \varphi) \right) \mid e^{{\rm i} \varphi} \in U(1) \right\}
\cong U(1).
\]

Following the same method as for the other phases, we conclude that
\[
\ell (\mathbf{w}, A) = \left\langle \! \left\langle \left[\widehat{\mathbf{w}}, A \right], \left[\widehat{\mathbf{w}}, A \right] \right\rangle \! \right\rangle = \mathbf{w}^ \mathsf{T} \mathbf{J}(A) \mathbf{w},
\]
where $\widehat{\mathbf{w}} = (\partial_z R) R ^{-1}$ and 
$\mathbf{J}(A)$ is given in \eqref{J_A}. Thus,
the Euler-Poincar\'e equations \eqref{EP_equations} become in this case
\begin{equation*} 
\partial_z \frac{\delta l}{\delta \mathbf{w}} +
\frac{\delta l}{\delta \mathbf{w}} \times \mathbf{w} = 
2\overrightarrow{\operatorname{Re}\left[\frac{\delta l}{\delta A}^*,  A\right]^{\phantom{4}}\!\! } ,\quad 
\partial _z A=  \left[\widehat{\mathbf{w}}, A \right],
\end{equation*}
and we have
\[ 
\frac{\delta l}{\delta \mathbf{w}}=2 \mathbf{J}(A) \mathbf{w}, \quad 
\frac{\delta l}{\delta A}= \left[\left[\widehat{\mathbf{w} },A\right] 
\Gamma ,\widehat{\mathbf{w} }\right] + 
\left[\left[\widehat{\mathbf{w} },A\right] 
\Gamma ,\widehat{\mathbf{w} }\right]^\mathsf{T}.
\]
As in the case of the $\Omega_6$-phase we obtain the following result.

\begin{theorem}
The three functions $h, j_m ,  \mathbf{J}_3$ form a completely 
integrable system on the six dimensional coadjoint orbit
$\Omega_8$.
\end{theorem}

\section{Acknowledgments}
This paper was finished when the authors 
were members of the programs ``Mathematics of Liquid Crystals" and
``Mathematical Modelling and Analysis of Complex Fluids and Active Media in Evolving Domains" at the Isaac Newton Institute for Mathematical 
Sciences, Cambridge, UK, January 7 - August 23, 2013. We express our 
thanks to the Institute and the organizers of the program for
the invitation and a very fruitful atmosphere conducive to 
work and collaboration.

{\footnotesize

\bibliographystyle{new}

\end{document}